\documentclass[10pt]{amsart}

\usepackage[english]{babel}
\usepackage[ansinew]{inputenc}
\usepackage{amsmath, amsthm, amssymb, graphicx, mathrsfs}
\usepackage{color}

\newcommand{\Zb}{\mathbb Z}

\makeatletter %Definition of \overrightharpoonup
   \def\rightharpoonupfill@{\arrowfill@\relbar\relbar\rightharpoonup}
   \newcommand{\overrightharpoonup}{%
   \mathpalette{\overarrow@\rightharpoonupfill@}}
\makeatother 
\makeatletter %Definition of \overleftharpoonup
   \def\leftharpoonupfill@{\arrowfill@\leftharpoonup\relbar\relbar}
   \newcommand{\overleftharpoonup}{%
   \mathpalette{\overarrow@\leftharpoonupfill@}}
\makeatother 
\def\epsilon{\varepsilon}
\def\phi{\varphi}
\def\Pr{\mathbb P}

\newcommand{\Curr}{\mbox{Curr}}
\newcommand{\Out}{\mbox{Out}}
\newcommand{\Aut}{\mbox{Aut}}

\newcommand{\PCN}{\Pr\Curr(\FN)}

\newcommand{\FN}{F_N}   % F ou F_n ou F_N ?
\newcommand{\cvn}{\mbox{cv}_N}
\newcommand{\cvnbar}{\overline{\mbox{cv}}_N}
\newcommand{\CVN}{\mbox{CV}_N}
\newcommand{\CVNbar}{\overline{\mbox{CV}}_N}
\newcommand{\R}{\mathbb R}
\newcommand{\Z}{\mathbb Z}

\usepackage{amsmath,amsthm,amssymb}

\newtheorem{theor}{Theorem}

\newtheorem{thm}{Theorem}[section]
\newtheorem{lem}[thm]{Lemma}

\newtheorem{prop}[thm]{Proposition}
\theoremstyle{definition}
\newtheorem{defn}[thm]{Definition}
\newtheorem{notation}[thm]{Notation}
\newtheorem{conv}[thm]{Convention}
\newtheorem{rem}[thm]{Remark}
\newtheorem{exmp}[thm]{Example}
\newtheorem{prop-defn}[thm]{Proposition-Definition}

\newtheorem{prob}[thm]{Problem}

\newtheorem*{claim*}{Claim}

\title{Spectral rigidity of automorphic orbits in free groups}

\author{Mathieu Carette}
\address{\tt SST/IRMP, 
Chemin du Cyclotron 2, bte L7.01.01, 
1348 Louvain-la-Neuve,
Belgium}
\email{\tt mathieu.carette@uclouvain.be}

\author{Stefano Francaviglia}

\address{\tt Dipartimento di Matematica of the University of Bologna, 
Pizza Porta S. Donato 5, 40126 Bologna, Italy}
\email{\tt francavi@dm.unibo.it }

\author{Ilya Kapovich}
\address{\tt Department of Mathematics, University of Illinois at
 Urbana-Champaign, 1409 West Green Street, Urbana, IL 61801, USA
 \newline http://www.math.uiuc.edu/\~{}kapovich/} \email{\tt
  kapovich@math.uiuc.edu}

\author{Armando Martino}
\address{\tt School of Mathematics,
University of Southampton;
Highfield, Southampton,
SO17 1BJ, United Kingdom }
\email{\tt A.Martino@soton.ac.uk}

\begin{document}
	
\begin{abstract}
It is well-known that a point $T\in \cvn$ in the (unprojectivized) Culler-Vogtmann Outer space $\cvn$ is uniquely determined by its \emph{translation length function} $||.||_T:F_N\to\mathbb R$.
A subset $S$ of a free group $F_N$ is called \emph{spectrally rigid} if, whenever $T,T'\in \cvn$ are such that $||g||_T=||g||_{T'}$ for every $g\in S$ then $T=T'$ in $\cvn$.
By contrast to the similar questions for the Teichm\"uller space, it is known that for $N\ge 2$ there does not exist a finite spectrally rigid subset of $F_N$. 

In this paper we prove that for $N\ge 3$  if $H\le \Aut(F_N)$ is a
subgroup that projects to a nontrivial normal subgroup in $\Out(F_N)$ then the $H$-orbit of an arbitrary nontrivial element $g\in F_N$ is spectrally rigid. We also establish a similar statement for $F_2=F(a,b)$,   provided that $g\in F_2$ is not conjugate to a power of $[a,b]$.
\end{abstract}

\thanks{The first author is a Postdoctoral Researcher of the F.R.S.-FNRS (Belgium). The third author was supported by the NSF
  grant DMS-0904200}


\maketitle

\tableofcontents

\section{Introduction}

The phenomenon of \emph{marked length spectrum rigidity} plays an important role in Riemannian geometry and adjacent areas. If $M$ is a closed manifold with a Riemannian metric $\rho$ of negative (but not necessarily constant) curvature, the associated \emph{marked length spectrum} is the function $\ell_\rho: G\to \mathbb R$, where $G=\pi_1(M)$ and where for $\gamma\in G$ $\ell_\rho(\gamma)$ is the shortest length with respect to $\rho$ among all free homotopy representatives of $\gamma$. It is easy to see that $\ell_\rho(\gamma)=\ell_\rho(\gamma_1^{-1}\gamma\gamma_1)$ for all $\gamma, \gamma_1\in G$. Thus $\ell_\rho$ may be also viewed as a function from the set of conjugacy classes in $G$ to $\mathbb R$.  One can also think of $\ell_\rho$ as the "translation length" function on $G$. Namely,  let $X=(\widetilde M, d_\rho)$, where $d_\rho$ is the distance function on $\widetilde M$ corresponding to the lift $\widetilde \rho$ of $\rho$ to $\widetilde M$.  Then the natural action of $G$ on $\widetilde M$ by covering transformation is an action of $G$ by isometries on $X$ (and thus can be thought of as a representation $G\to Isom(X))$ and for every $\gamma\in G$ we have
\[
\ell_\rho(\gamma)=\inf_{x\in X} d_\rho(x,\gamma x)=\min_{x\in X} d_\rho(x,\gamma x)
\]
is the \emph{translation length} of $\gamma$ as the isometry of $X$.

The \emph{Marked Length Spectrum Rigidity Conjecture} (MLSRC) states
that knowing the function $\ell_\rho$ uniquely determines the isometry
type of $(M,\rho)$. More precisely, the conjecture says that if $\rho,\rho'$ are two smooth
negatively curved Riemannian metrics on $M$ such that
$\ell_\rho=\ell{\rho'}$ then there exists an isometry from $(M,\rho)$
to $(M,\rho')$ which is isotopic to the identity. This conjecture is
known to hold in some special cases (more on this below), but is still
open in full generality.

There are also various generalizations of MLSRC to other contexts,
such as allowing more general types of metrics on $M$. There are also
generalizations  with the set-up where,  given two representations
$\tau: G\to Isom(X_1)$ and $\tau_2:G\to Isom(X_2)$ with the same marked
length spectrum $\ell_{\tau_1}=\ell_{\tau_2}: G\to\mathbb R$ (where
$X_1$ and $X_2$ are required to satisfy various kinds of negative or
non-positve curvature conditions). Then the desired conclusion of
MLSRC is that there  exists an isometry $X_1\to X_2$ conjugating
$\rho_1$ to $\rho_2$.  MLSRC is known to hold for surfaces, including
various generalizations of the types of metrics $\rho$ included under
consideration, see~\cite{Cr90,Otal,CFF}. There are also a number of
known results establishing versions of MLSRC  for representations
$\rho_1: G\to Isom(X_1)$ , $\rho_2:G\to Isom(X_2)$ where $X_1,X_2$ are
allowed to be higher dimensional, but there are more significant
restrictions on the geometry of $X_1,X_2$ than in the results about
MLSRC for surfaces (see, for example, \cite{CEK,HP,Kim01a,Kim01b,DK,Kim04}).  However, the original version of MLSRC is still mostly open (except for rather special classes of metrics) in dimensions bigger than two. We refer the reader to the survey~\cite{Cr04} for a more extended discussion on the topic.

In any context where MLSRC is known to hold, it is natural to ask if there are smaller subsets of $G$ such that knowing the restriction of the marked length spectrum to such a subset uniquely determines the entire marked length spectrum. Namely, for a given class of length functions $\ell:G\to\mathbb R$  where MLSRC is known to hold, say that a subset $S\subseteq G$ is \emph{spectrally rigid} if whenever $\ell, \ell'$ are two length functions  from the class in question such that $\ell|_S=\ell'|_S$ then $\ell=\ell'$. 
For closed surfaces with metrics of constant curvature $-1$ the situation is particularly well-behaved.
Thus it is known (see, for example,~\cite{FLP}) that if $\Sigma$ is a
closed oriented surface of genus $\ge 2$, then there exist elements
$h_1, \dots, h_{6g-5}\in G=\pi_1(\Sigma)$ such that whenever
$\rho_1,\rho_2$ are two points in the Teichmuller space $\mathcal
T(\Sigma)$ (i.e. marked hyperbolic metrics on $\Sigma$) such that
$\ell_{\rho_1}(h_i)=\ell_{\rho_2}(h_i)$ for  $i=1,\dots , 6g-5$ then
$\ell_{\rho_1}=\ell_{\rho_2}$ and $\rho_1=\rho_2$ in  $\mathcal
T(\Sigma)$. Thus the subset $\{h_1,\dots, h_{6g-5}\}\subseteq G$
``spectrally rigid" with respect to the class of hyperbolic metrics on
$\Sigma$.  We believe that looking for ``small'' spectrally rigid sets
in other situations, where MLSRC is known to hold,  is an interesting
general problem worthy of further study.

If $G$ is a finitely generated group acting by isometries on an $\mathbb R$-tree $X$, there is also a naturally associated \emph{translation length function} $||.||_X=\ell_X:G\to \mathbb R$, where for $g\in G$
\[
\ell_X(g)=\inf_{x\in X} d_X(x,g x)=\min_{x\in X} d_X(x, gx).
\]
It is well-known~\cite{CM,Pau89,Chis} that under some mild extra assumptions (which are satisfied, in particular, for the Outer space context discussed below), MLSRC holds, that is, knowing the function $\ell_X$ uniquely determines $X$ and the action of $G$ on $X$ (up to a $G$-equivariant isometry). 

For a free group $F_N$ (where $N\ge 2$) the \emph{Culler-Vogtmann
  Outer space} $\cvn$ is an analog of the Teichmuller space of a
hyperbolic surface. The space $\cvn$ consists of minimal free discrete
isometric actions of $F_N$ on $\mathbb R$-trees, considered up to
$F_N$-equivariant isometry. Every element $T\in \cvn$ arises as the
universal cover of a finite graph $\Gamma$, whose fundamental group is
identified with $F_N$ via a particular isomorphism, where edges of
$\Gamma$ are given positive real lengths and their lifts to $T$ are
given the same lengths. There is an important subset, $\CVN\subseteq
\cvn$, consisting of all $T\in \cvn$ such that the quotient metric
graph $T/F_N$ has volume 1. The space $\CVN$ is the
\emph{projectivized Culler-Vogtmann Outer space}, which was introduced
by Culler and Vogtmann \cite{CV}, before the introduction of $\cvn$ . Both $\cvn$ and $\CVN$ play an important role in the study of $Out(F_N)$.
We say that a subset $R\subseteq F_N$ is \emph{spectrally rigid} if
whenever $T_1, T_2\in \cvn$ are such that $||g||_{T_1}=||g||_{T_2}$
for every $g\in R$ then $T_1=T_2$ in $\cvn$. As noted above, $R=F_N$
is spectrally rigid. Moreover, for every $T\in \cvn$ the translation
length function $||.||$ is constant on every conjugacy class in
$F_N$. Thus if $R$ is chosen to consist of representatives of all
conjugacy classes in $F_N$, then $R$ is spectrally rigid.

A surprising result of Smillie and Vogtmann~\cite{SV} shows that there
does not exist a finite spectrally rigid subset of $F_N$, where $N\ge
3$. A result of Cohen, Lustig and Steiner~\cite{CLS} establishes the same fact for $N=2$. In particular, it is proved in \cite{SV}  that for any finite subset $R\subseteq
F_N$, where $N\ge 3$, there exists a one-parametric family $(T_t)_{t\in [0,1]}$ of
distinct points of $\CVN$ such that for every $t\in [0,1]$ the length
functions $||.||_{T_0}$ and $||.||_{T_t}$ agree on $R$.  A similar
statement follows from the recent work of  Duchin, Leininger, and Rafi
for flat metrics on surfaces of finite type~\cite{DLR}. Moreover, the paper \cite{DLR} gives a complete characterization of when a set of simple closed curves on a finite type surface is spectrally rigid with respect to the space of flat metrics on that surface.

In view of the results of~\cite{CLS,SV} it becomes interesting to look
for infinite but "sparse" spectrally rigid subsets of $F_N$. The study
of this topic was initiated by Kapovich in \cite{Ka4}. Namely, it is
proved in~\cite{Ka4} that if $A=\{a_1,\dots, a_N\}$ is a free basis of
$F_N$ then  almost every trajectory of the simple non-backtracking
random walk on $F_N$ with respect to $A$ yields a spectrally
rigid subset of $F_N$. Recently Brian
Ray~\cite{Ray} proved that for any $\phi\in \Aut(F_N)$ (where $N\ge
1$) and for any $g\in F_N$ the set $\langle \phi \rangle
g=\{\phi^n(g): n\in \mathbb Z\}$ is not spectrally rigid in
$F_N$. Ongoing work of Ray (unpublished) also shows that for any
finite collection $H_1,\dots, H_m\le F_N$ of finitely generated
subgroups of infinite index in $F_N$ (where $N\ge 2$) the set
$H_1\cup\dots \cup H_m$ is not spectrally rigid in $F_N$.

In the present paper we obtain a very different class of examples of
spectrally rigid subsets of free groups from those constructed in
\cite{Ka4}. We say that a subgroup $H\le \Aut(F_N)$ (where $N\ge 2$)
is \emph{ample} if the image of $H$ in $\Out(F_N)$ contains
an infinite normal subgroup of $\Out(F_N)$. It is
well-known~\cite{Z,Cu} that every finite subgroup of $\Out(F_N)$ comes
from a group of simplicial automorphisms of a finite connected graph
without degree-one vertices and with fundamental group $F_N$. From
here it is not hard to show that for $N\ge 3$ every nontrivial normal
subgroup of $\Out(F_N)$ is infinite. For the case $N=2$ the group
$\Out(F_2)$ possesses a unique nontrivial finite normal subgroup,
namely the center $Z(\Out(F_N))$ which is a cyclic group of order $2$
generated by the ``hyper-elliptic involution'' of $F_2=F(a,b)$,
$a\mapsto a^{-1}$, $b\mapsto b^{-1}$.

Our main result is:

\begin{theor}\label{thm:A}
Let $N\ge 2$ and let $H\le \Aut(F_N)$ be an ample subgroup.
Let $g\in F_N$ be an arbitrary nontrivial element; in the case $N=2$
we also assume that $g\in F_2=F(a,b)$ is not conjugate to a nonzero power of $[a,b]$ in $F_2$.

Then the orbit $Hg=\{\phi(g) : \phi\in H\}$ is a spectrally rigid subset of $F_N$.
\end{theor}
Theorem~\ref{thm:A} applies, for example, to the cases where
$H=\Aut(F_N)$ (with $N\ge 2$) or where $H\le \Aut(F_N)$ is the \emph{Torelli subgroup}
(with $N\ge 3$), that is, $H$ is the set of all elements of $\Aut(F_N)$ that act as the identity map on the abelianization  $\Z^N$ of $F_N$.
Theorem~\ref{thm:A} also immediately implies that for $N\ge 3$ any
$\Aut(F_N)$-invariant subset of $F_N$ with more than one element is
spectrally rigid in $F_N$.

For $F_2=F(a,b)$ it is well-known that for every $\phi\in \Aut(F_2)$
the element $\phi([a,b])$ is conjugate to $[a,b]^{\pm 1}$ in $F_2$,
which easily implies that $\Aut(F_2)[a,b]$ is not spectrally rigid in
$F_2$. However, as Theorem~\ref{thm:A} shows, this example is
essentially the only obstruction for spectral rigidity of an
$\Aut(F_N)$-orbit of a nontrivial element in the case $N=2$.

As noted above, for any $T\in\cvn$ $||.||_T$ is a class function and
thus it can be viewed as a function on the set $\mathcal C_N$ of all
the conjugacy classes of elements of $F_N$. The notion of a spectrally
rigid set also naturally translates to subsets of $\mathcal
C_N$. However, the results of Brian Ray about finitely generated
subgroups of $F_N$ mentioned above, and the results of \cite{Ka4}
about random walk trajectories are more naturally formulated in the
setting of subsets of $F_N$. For that reason we adhere to that setting
in this paper.

The first step in the proof of Theorem~\ref{thm:A} is to establish Theorem~\ref{thm:prim} which says
that the set $\mathcal P_N$ of all the primitive elements in $F_N$ is
a spectrally rigid subset in $F_N$ for every $N\ge 2$. Recall that an
element $g\in F_N$ is \emph{primitive} if it belongs to some free
basis of $F_N$. Thus $\mathcal P_N=\Aut(F_N)g$ where $g\in F_N$ is any
primitive element.  Theorem~\ref{thm:prim} is derived from the results
of Francaviglia and Martino~\cite{FraMa} about extremal Lipschitz
distortions between two arbitrary points in $\cvn$.   A key fact there
is that for any $T,T'\in \cvn$ the ``extremal Lipschitz distortion''
$D(T,T'):=\sup_{g\in F_N, g\ne 1} \frac{||g||_{T'}}{||g||_T}$ is
actually a maximum which is realized by an element $g$ from some
finite subset $\mathbf U_T\subseteq F_N$ depending only on $T$. Moreover,
the explicit description of elements of $\mathbf U_T$ in~\cite{FraMa}
shows that they are all primitive, so that $\mathbf U_T\subseteq
\mathcal P_N$. From here it is easy to see that if
$||g||_T=||g||_{T'}$ for every $g\in \mathcal P_N$ then
$D(T,T')=D(T',T)=1$ and hence $T=T'$ in $\cvn$. 
A more careful version of the above argument yields the following "relative rigidity" result:

\begin{theor}\label{thm:C}
Let $T\in \cvn$ be arbitrary. There exists a finite set $S$ (depending on $T$) of primitive elements in $F_N$ with the following property:
Whenever $T'\in \cvn$ is such that $||g||_{T'}=||g||_T$ for every $g\in S$ then $T=T'$ in $\cvn$.
\end{theor}
In fact, the proof of Theorem~\ref{thm:C} shows that we can take $S=\mathbf U_T$.

After Theorem~\ref{thm:prim} is established, we derive Theorem~\ref{thm:A} from Theorem~\ref{thm:prim} using the machinery of \emph{geodesic currents} on free groups, and particularly exploiting the \emph{geometric intersection form} between trees and currents, constructed in \cite{Ka2,KL2}.
A geodesic current is a measure-theoretic analog of the notion of the
conjugacy class in a (word-hyperbolic) group. Geodesic currents were
introduced and studied by Bonahon~\cite{Bo86,Bo88} in the context of
hyperbolic surfaces and the Teichmuller space, where they turned out
to be quite useful. A \emph{geodesic current} on a free group $F_N$ is
a positive Radon measure $\mu$ on $\partial^2 F_N=\partial F_N \times
\partial F_N- diag$ which is "flip"-invariant and $F_N$-invariant. The
space $\Curr(F_N)$ of all geodesic currents on $F_N$ comes equipped
with a natural weak-* topology making it into a locally compact space,
and with a natural left $\Out(F_N)$-action by linear transformations.
The theory of geodesic currents on free groups has been actively
developed in the last several years by Kapovich~\cite{Ka1,Ka2,Ka3,Ka4} and Kapovich-Lustig~\cite{KL1,KL2,KL3,KL4}
(see also ~\cite{BF08,CP,Ha,Fra,KN1,KN2,KN3} for other recent applications of currents). The space $\Curr(F_N)$ turns out to be a natural counterpart for the Outer space $\cvn$, and, more generally, the closure $\cvnbar$ of $\cvn$.  
The closure $\cvnbar$ of $\cvn$ (with respect to equivariant Gromov-Hausdorff convergence topology) is known to consist of all the minimal \emph{very small} isometric actions of $F_N$ on $\R$-trees. The Outer space $\cvn$ is an open $\Out(F_N)$-invariant dense subset of $\cvnbar$. It is again well-known that any point $T\in\cvnbar$ is uniquely determined by its translation length function $||.||_T: F_N\to \R$.

The interaction between $\cvnbar$ and $\Curr(F_N)$ is given by the \emph{geometric intersection form}
\[
\langle \cdot\, , \ \cdot \rangle: \cvnbar\times \Curr(F_N)\to \mathbb R_{\ge 0},
\] 
constructed in~\cite{KL2}. The intersection form has a number of useful properties, such as being continuous, $Out(F_N)$-invariant, $\mathbb R_{\ge 0}$-linear with respect to the second argument and $\mathbb R_{\ge 0}$-linear with respect to the first argument. Another key property of the intersection form, relating it to marked length spectra, is that for every $T\in \cvnbar$ and every $g\in F_N, g\ne 1$ we have
\[
\langle T, \eta_g\rangle=||g||_T
\]
where $\eta_g\in \Curr(F_N)$ is the so-called "counting" current
defined by $g$. This last property is crucial in deriving
Theorem~\ref{thm:A} from
Theorem~\ref{thm:prim}. 

As noted above, the first step in the proof of Theorem~\ref{thm:A} is
establishing Theorem~\ref{thm:prim}, which is accomplied via analyzing
extremal Lipschitz distortion between arbitrary points in
$\cvn$. After that the proof of Theorem~\ref{thm:A} (for $N\ge 3$)
proceeds as follows. It is known~\cite{KL1} that for $N\ge 3$ there
exists a unique minimal closed nonempty $\Out(F_N)$ subset of the
space $\PCN$ of projectivized geodesic currents: this subset $\mathbb
M_N$, called the \emph{minimal set}, is exactly the closure in $\PCN$
of the set of projective classes $[\eta_a]$ for all the primitive
elements $a\in F_N$. Now let $H\le \Aut(F_N)$ be an ample
subgroup and let $g\in F_N$ be a nontrivial element. First, using a
powerful recent result of Handel and Mosher~\cite{HM} we conclude that
$H$ contains a ``fully irreducible'' (also known as ``iwip'', which
stands for ``irreducible with irreducible powers'') element
$\psi\in H$. Then, using the assumption about $H$ being ample, and
the dynamics of the action of fully irreducible automorphisms on $\PCN$, we
conclude that the closure in $\PCN$ of the orbit $H[\eta_g]$ contains
the minimal set $\mathbb M_N$. Then, using the intersection form
$\langle \cdot\, , \ \cdot \rangle$ mentioned above, we conclude that
for $T\in \cvn$ knowing all the translation lengths
$||\phi(g)||_T=\langle T, \eta_{\phi(g)}\rangle$, where $\phi\in H$,
allows us to recover $||a||_T=\langle T, \eta_a\rangle$ for all
primitive elements $a\in F_N$. Finally, by applying
Theorem~\ref{thm:prim} we conclude that the set $\{\phi(g): \phi\in
H\}$ is spectrally rigid in $F_N$.

The proof in \cite{Ka4} that a random trajectory of the simple non-backtracking random walk on $F_N$ yields a spectrally rigid subset of $F_N$, was also based on using geodesic currents and the geometric intersection form. A key fact in that proof was that the $\R_{\ge 0}$-linear span of the set of counting currents $\eta_{w_n}$, where $(w_n)_{n\ge 1}$ was a random trajectory of the walk, is dense in $\Curr(F_N)$. A similar line of reasoning cannot be used to prove Theorem~\ref{thm:A}. Indeed, if $A$ is a free basis of $F_N$ and $v$ is a freely reduced word over $A$, such that every freely reduced word of length $2$ over $A$ occurs in $v$ as a subword, then $v$ cannot be a subword of a cyclic word representing a primitive element in $F_N$. This implies that $v$ has "weight $0$" (see~\cite{Ka2} for the relevant terminology) in any finite linear combination of counting currents of primitive elements. Hence, by continuity,  $v$ has "weight $0$" in every current from the closure $Z$ of the linear span of the counting currents of all primitive elements, and hence $Z$ is a proper subset of $\Curr(F_N)$.

In Section~\ref{sect:op} we discuss a number of open problems motivated by the results of this paper, including questions about spectral rigidity properties of subgroups of $\Out(F_N)$, questions about spectral rigidity with respect to the closure $\cvnbar$ of $\cvn$, and other problems.

We are grateful to Martin Lustig for useful conversations and to Jing Tao for providing Example~\ref{exmp:Tao}.  We also thank the referee for useful suggestions.
	
	\section{Preliminaries}

\subsection{Graphs and graph-related conventions}

\begin{conv}[Graphs]

A \emph{graph} is a 1-complex. The set of $0$-cells of a graph
$\Delta$ is denoted $V\Delta$ and its elements are called
\emph{vertices} of $\Delta$. The closed 1-cells of a graph $\Delta$ are
called \emph{topological edges} of $\Delta$. The set of all
topological edges is denoted $E_{top}\Delta$. 

The interior of every topological edge is
homeomorphic to the interval $(0,1)\subseteq \mathbb R$ and thus
admits exactly 2 orientations (when considered as a 1-manifold). We
call a topological edge endowed with the choice of an orientation on
its interior an \emph{oriented edge} of $\Delta$. The set of all oriented edges
of $\Delta$ is denoted $E\Delta$. For an oriented edge $e\in E\Delta$
changing its orientation to the opposite produces another oriented
edge of $\Delta$ denoted $e^{-1}$ and called the \emph{inverse} of
$e$. Thus ${}^{-1}:E\Delta\to E\Delta$
is a fixed-point-free involution.

For every oriented edge $e$ of $\Delta$ there are naturally defined
(and not necessarily distinct)
vertices $o(e)\in V\Delta$, called the \emph{origin} of $e$, and
$t(e)\in E\Delta$, called the \emph{terminus} of $e$, satisfying
$o(e^{-1})=t(e)$, $t(e^{-1})=o(e)$. 
\end{conv}

An \emph{orientation} on a graph $\Delta$ is a partition $E\Delta=E^+\Delta\sqcup E^-\Delta$, where for every $e\in E\Delta$ one
of the edges $e, e^{-1}$ belongs to $E^+\Delta$ and the other edge belongs to  $E^-\Delta$.

\begin{defn}[Paths]\label{defn:path}
A \emph{simplicial path} or \emph{edge-path} $\gamma$ of \emph{simplicial length} $n\ge 1$
in $\Delta$ is a sequence of oriented edges 
\[
\gamma=e_1,\dots, e_n
\]
such that $t(e_i)=o(e_{i+1})$ for $i=1,\dots, n-1$. We say that
$o(\gamma):=o(e_1)$ is the \emph{origin} of $\gamma$ and that
$t(\gamma)=t(e_n)$ is the \emph{terminus} of $\gamma$. We also regard
a $\gamma=v\in V\Delta$ as an simplicial path in $\Gamma$ of
simplicial length $0$ with $o(\gamma)=t(\gamma)=v$. A simplicial path is
called \emph{reduced} if it does not contain a back-tracking, that is
a path of the form $ee^{-1}$, where $e\in E\Delta$.

\end{defn}

\subsection{Outer space}

The Culler-Vogtmann Outer space, introduced by Culler and Vogtmann in
a seminal paper~\cite{CV}, is a free group analogue of the
Teichm\"uller space of a closed surface of negative Euler
characteristic. We briefly review some basics definitions and facts
about the Outer space, and refer the reader to \cite{CV,BF93,BV,Gui1,V08} for more
detailed background information on the topic.

\begin{defn}[Non-projectivized Outer Space]
Let $F_N$ be a finitely generated free group of rank $N\ge 2$.

The \emph{non-projectivized outer space} $\cvn$ consists of all
minimal free and discrete isometric actions of $F_N$ on $\mathbb
R$-trees. Two such trees $T_1,T_2$ are considered equal in $\cvn$ if
there exists an $F_N$-equivariant isometry between them.
The space $\cvn$ is endowed with the equivariant Gromov-Hausdorff
convergence topology.
\end{defn}

For $T\in \cvn$ and $c>0$ denote by $cT\in \cvn$ the tree that
coincides with $T$ as a topological space and has the same $F$-action,
but where the metric is multiplied by $c$.

A basic fact in the theory of Outer space states that every $T\in \cvn$ is uniquely determined by its
\emph{translation length function} $||.||_T: F_N\to \mathbb R$, where for
every
$g\in F_N$
\[
||g||_T=\min_{x\in T} d_T(x,gx)
\]
is the \emph{translation length} of $g$.

\begin{prop}\label{prop:mlsr}
Let $T_1,T_2\in \cvn$ be such that $||.||_{T_1}=||.||_{T_2}$, that is
$||g||_{T_1}=||g||_{T_2}$ for every $g\in F_N$. Then there exists an
$F_N$-equivariant isometry between $T_1$ and $T_2$, so that $T_1=T_2$
in $\cvn$.
\end{prop}

Proposition~\ref{prop:mlsr} is a special case of a much more general
fact. Thus it is known~\cite{Pau89} that for a finitely generated
group $G$ any minimal irreducible (without a global fixed end) isometric action of
$G$ on an $\mathbb R$-tree $T$ is uniquely determined (up to a
$G$-equivariant isometry) by the translation length function $||.||_T:
G\to\mathbb R$ for this action.

Note that
$||g||_T=||hgh^{-1}||_T$ for every $g,h\in F_N$. Thus $||.||_T$ can
be thought of as a function on the set of conjugacy classes in $F_N$.
The space $\cvn$ comes equipped with a natural right $Out(F_N)$-action
by homeomorphisms. At the length function level, if $\phi\in Out(F_N)$,
$T\in \cvn$ and $g\in F_N$ we have
\[
||g||_{T\phi}=||\phi(g)||_T.
\]
It is known that the equivariant Gromov-Hausdorff topology on $\cvn$
coincides with the pointwise convergence topology at the level of
length functions. Thus for $T_n,T\in \cvn$ we have $\lim_{n\to\infty}
T_n=T$ if and only if for every $g\in F$ we have $\lim_{n\to\infty} ||g||_{T_n}=||g||_T$.

\begin{defn}[Projectivized Outer Space]
Denote by $\CVN$ the subset of $\cvn$ consisting of all $T\in \cvn$
such that the quotient graph $T/F_N$ has volume $1$.

The space $\CVN$ is a closed $Out(F_N)$-invariant subset of $\cvn$ and
it is called the \emph{projectivized Outer Space} of $F_N$.
\end{defn}

It is known that $\CVN$ is $Out(F_N)$-equivariantly homeomorphic to $\cvn/\sim$,
where $T_1\sim T_2$ if there is $c>0$ such that $T_2=cT_1$ in $\cvn$. This fact
justifies the term ``projectivized Outer Space''. For $T\in \cvn$ we denote by $[T]$ the $\sim$-equivalence class of $T$ and call $[T]$  \emph{the projective class of $T$}.

Points of $\cvn$ have a more explicit combinatorial description as
``marked metric graph structures'' on $F$:

\begin{defn}[Metric graph]
A  \emph{metric graph} is a graph $\Delta$ endowed with a \emph{metric
  structure} $\mathcal L$, that is, a function $\mathcal L:E\Delta\to
(0,\infty)$ such that for every $e\in E\Delta$ we have $\mathcal
L(e)=\mathcal L(e^{-1})$. The number $\mathcal L(e)$ is called the
\emph{length} of $e$ with respect to $\mathcal L$. 

For a metric graph $(\Delta,\mathcal L)$ its \emph{volume} is defined 
as
\[
vol_\mathcal L(\Delta)=\frac{1}{2}\sum_{e\in E\Delta} \mathcal L(e)=\sum_{e\in E^+\Delta} \mathcal L(e),
\]
where $E\Delta=E^+\Delta\sqcup E^-\Delta$ is any orientation on $\Delta$.

\end{defn}

\begin{defn}[Marking]
Let $F_N$ be a free group of finite rank $N\ge 2$. A \emph{marking} or
a \emph{marked graph structure} on
$F_N$ is an isomorphism $\alpha:F_N\to \pi_1(\Gamma)$ where $\Gamma$
is a finite connected graph with the first Betti number equal to $N$
and such that $\Gamma$ has no degree-1 and no degree-2 vertices. 
\end{defn}

\begin{defn}[Marked metric graph]
A \emph{marked metric graph} or a \emph{marked metric graph structure}
on $F_N$ consists of a marking $\alpha:F_N\to \pi_1(\Gamma)$ on $F_N$
together with a metric graph structure $\mathcal L$ on $\Gamma$.
\end{defn}

\begin{conv}\label{conv:mgs}
Let $(\alpha, \mathcal L)$ be a marked metric graph structure on
$F_N$, where $\alpha: F_N\to \pi_1(\Gamma,p)$ is a marking and where
$\mathcal L$ is a metric structure on $\Gamma$.

Then $(\alpha, \mathcal L)$
defines a point $T\in \cvn$ as follows. Topologically, let $T=\widetilde
\Gamma$, with an action of $F_N$ on $T$ via $\alpha$. We lift the metric
structure $\mathcal L$ from $\Gamma$ to $T$ by giving every edge in
$T$ the same length as that of its projection in $\Gamma$. This makes
$T$ into an $\mathbb R$-tree equipped with a minimal free and discrete
isometric action of $F_N$. (The assumption that $\Gamma$ has no
degree-1 vertices guarantees that the action of $F_N$ on $T$ is minimal). Thus $T\in \cvn$ and in this situation we
will sometimes use the  notation $T=(\alpha, \mathcal L)\in \cvn$.
Note that $T/F=\Gamma$. 

Moreover, it is not hard to see that every
point of $\cvn$ arises in this fashion and that $\cvn$ is exactly
the set of all those $T=(\alpha, \mathcal L)\in \cvn$ where $(\alpha,
\mathcal L)$ is a marked metric graph structure on $F_N$ with
$vol_\mathcal L(\Gamma)=1$.

For this reason we will also think of elements $T\in \cvn$ as metric
graphs, and the default assumption will be that every vertex of $T$
has degree $\ge 3$.
\end{conv}

The following useful proposition is an immediate corollary of the definitions:

\begin{prop}\label{circuit}
Let $T\in \cvn$ be realized by a marked metric graph structure 
$(\alpha: F_N\to \pi_1(\Gamma), \mathcal L)$ on $F_N$, so that
$T=(\widetilde \Gamma, d_\mathcal L)$. Let $g\in
F_N, g\ne 1$. Let $\gamma_g$ be the unique immersed circuit in
$\Gamma$ obtained by reducing and cyclically reducing the edge-path
$\alpha(g)$ in $\Gamma$. 

Then $||g||_T$ is equal to the $\mathcal L$-length of $\gamma_g$.
\end{prop}

The outer space $\cvn$ has a natural closure $\cvnbar$ with respect to
the equivariant Gromov-Hausdorff convergence topology (or,
equivalently, with respect to the length function topology). It is
known~\cite{BF93,CL,Gui} that $\cvnbar$ consists precisely of (the
$F_N$-equivariant isometry classes of) all the \emph{very small}
minimal isometric actions of $F_N$ on $\mathbb R$-trees. The
$\Out(F_N)$-action on $\cvn$ naturally extends to an action on
$\cvnbar$ by homeomorphisms. Moreover, the projectivization $\CVNbar$
of $\cvnbar$ is compact and contains (a copy of) $\CVN$ as an open
dense $\Out(F_N)$-invariant subset. The space $\CVNbar$ is sometimes
called the \emph{Thurston compactification} of $\CVN$.

\section{Extremal Lipschitz distortions and rigidity of the set of primitive elements}

For $N\ge 2$ we denote by $\mathcal P_N$ the set of all primitive elements in $F_N$.

\begin{notation}\label{not:D}
Let $T\in \cvn$ and $T'\in\cvnbar$. 
Denote 
\[
D(T,T'):=\sup_{g\in F_N, g\ne 1} \frac{||g||_{T'}}{||g||_T}.\tag{\ddag}
\]
\end{notation}	

\begin{notation}[Almost simple curves]
Let $T\in \cvn$ and let $\Gamma=T/F_N$ be the quotient metric graph, with the metric structure $\mathcal L$ coming from $T$. Thus $F_N$ is naturally identified with $\pi_1(\Gamma)$  via an isomorphism $\alpha:F_N\to \pi_1(\Gamma)$ and $T=(\alpha,\mathcal L)$ in $\cvn$.

Let $\mathbf U_T\subseteq F_N$ be the set of all elements of $F_N$ corresponding (under $\alpha$) to the closed curves $\gamma$ in $\Gamma$ of one of the following types:
\begin{enumerate}
\item $\gamma$ is a nontrivial simple closed circuit in $\Gamma$;
\item $\gamma=\gamma_1\gamma_2$ is a concatenation of two nontrivial simple closed circuits $\gamma_1$, $\gamma_2$, each beginning and ending at a common vertex $v$, and such that $\gamma_1$, $\gamma_2$ do not  contain any common topological edges.  We refer to such $\gamma$ as \emph{figure-eight curves} in $T/F_N$.
\item $\gamma$ is a "barbell" circuit, that is $\gamma=\gamma_1 \beta \gamma_2 \beta^{-1}$ where $\gamma_i$ is a nontrivial simple closed circuit at a vertex $v_i$ of $\Gamma$ with $v_1\ne v_2$, where $\beta$ is a simple edge-path from $v_1$ to $v_2$ in $\Gamma$ and where $\gamma_1,\beta, \beta_2$ do not have any common topological edges. We refer to such $\gamma$ as \emph{barbell curves} in $T/F_N$. 
\end{enumerate} 

Note that, by construction, every element of $\mathbf U_T$ is primitive in $F_N$ and the set $\mathbf U_T$ is finite. Moreover $\# \mathbf U_T\le K(N)$ for some constant $K(N)$ depending only on $N$.  We call elements of $\mathbf U_T$ \emph{almost simple curves} for $T/F_N$.
\end{notation}	

We need the following fact established by Francaviglia and Martino in~\cite{FraMa}:
\begin{prop}\label{prop:FM}
Let $T, T'\in \cvn$ be arbitrary.
Then
\[
D(T,T')=\max_{g\in \mathbf U_T} \frac{||g||_{T'}}{||g||_T}.
\]

\end{prop}
Thus we see that the supremum in the definition of $D(T,T')$ in $(\ddag)$ is a maximum and it is achieved on one of elements from the finite subset $\mathbf U_T\subseteq \mathcal P_N$.

Proposition~\ref{prop:FM} quickly implies that the set of all primitive elements is spectrally rigid:

\begin{thm}\label{thm:prim}
Let $N\ge 2$. Then the set $\mathcal P$ of all primitive elements in
$F_N$ is spectrally rigid in $F_N$.
\end{thm}
\begin{proof}
Let $T, T'\in\cvn$ be such that $||g||_T=||g||_{T'}$ for every $g\in \mathcal P_N$.

Then
\[
D(T,T')=\sup_{g\in F_N, g\ne 1} \frac{||g||_{T'}}{||g||_T}=\max_{g\in \mathbf U_T} \frac{||g||_{T'}}{||g||_T} =1
\]
since  $\mathbf U_T\subseteq \mathcal P_N$.
Similarly,
\[
D(T',T)=\sup_{g\in F_N, g\ne 1} \frac{||g||_{T}}{||g||_{T'}}=\max_{g\in \mathbf U_{T'}} \frac{||g||_{T}}{||g||_{T'}} =1
\]
so that
\[
\inf_{g\in F_N, g\ne 1} \frac{||g||_{T'}}{||g||_T}=\sup_{g\in F_N, g\ne 1} \frac{||g||_{T'}}{||g||_T}=1.
\]
Thus $||g||_T=||g||_{T'}$ for every $g\in F_N$ and hence $T=T'$ in $\cvn$, as required.

\end{proof}

We now state a more precise (compared to Proposition~\ref{prop:FM}) statement summarizing some of the results of \cite{FraMa}:
\begin{prop}\label{prop:FM1}
Let $T,T'\in\cvn$ be arbitrary. Let $L:=D(T,T')$.
Then there exists an $F_N$-equivariant $L$-Lipschitz map $f:T\to T'$ with the following properties:
\begin{enumerate}
\item On each edge $e$ of $T$ the map $f$ is a linear map with constant stretch $D_e\ge 0$ (note that, because of $F_N$-equivariance, $D_{e_1}=D_{e_2}$  whenever $e_1$ and $e_2$ are in the same $F_N$-orbit of edges of $T$).
\item There exists $h\in \mathbf U_T$ such that for every edge $e$ in the axis $A_h$ of $h$ in $T$ we have $D_e=L$ and the restriction of $f$ to $A_h$ is injective (so that $f|_{A_h}$ is an $L$-homothety). 
\end{enumerate}
\end{prop}

We can now establish Theorem~\ref{thm:C} from the Introduction:
\begin{thm}\label{thm:rel}

Let $T, T'\in \cvn$ be such that $||g||_T=||g||_{T'}$ for every $g\in \mathbf U_T$. Then $T'=T$ in $\cvn$.
\end{thm}
\begin{proof}
Let $T,T'\in  \cvn$ be such that $||g||_T=||g||_{T'}$ for every $g\in \mathbf U_T$.
Hence by Proposition~\ref{prop:FM} $D(T,T')=1$. 
Let $f:T\to T'$ be an $F_N$-equivariant $1$-Lipschitz map provided by Proposition~\ref{prop:FM1}.  
In particular, $D_e\le 1$ for every edge $e$ of $T$.

We claim that $D_e=1$ for every edge $e$ of $T$. Indeed, suppose not, so that there exists an edge $e_0$ of $T$ with $D_{e_0}<1$. Let $y_0$ be the edge of $T/F_N$ which is the projection of $e_0$ to $T/F_N$. Then there exists  an immersed circuit $\gamma_0$ in $T/F_N$ passing through the edge $y_0$ such that $\gamma_0$ is an almost simple curve in $F_N/T$.  Let $g_0\in F_N$ correspond to $\gamma_0$, so that $g_0\in \mathbf U_T$. The fact that $D_{e_0}<1$ and that $D_e\le 1$ for every edge $e$ of $T$ implies that $||g_0||_{T'}<||g_0||_T$.  Again, this contradicts our assumptions on $T, T'$ and the fact that $g_0\in \mathbf U_T$.

We now claim that $f$ is injective.
Indeed, suppose not. Then there exist two distinct oriented edges $e_1, e_2$ of $T$ with a common initial vertex $v$ such that $f$ "folds" a non-degenerate initial segment of $e_1$ and a non-degenerate initial segment of $e_2$.  Let $y_1,y_2$ be the edges of $T/F_N$ which are the projections of $e_1$ and $e_2$ respectively to $T/F_N$. From the definition of an almost simple curve for $T/F_N$ it follows that there exists an immersed circuit $\gamma$ in $T/F_N$ containing $y_1^{-1} y_2$ as a subpath such that $\gamma$ is an almost simple curve in $T/F_N$. Thus $\gamma$ corresponds to $g\in \mathbf U_T$. By $F_N$-equivariance of $f$ we may assume that $e_1^{-1}e_2$ is a subpath of the axis of $g$ in $T$.  Since $D_e=1$ for every edge $e$ of $T$ and since $f$ folds nondegenerate initial segments of $e_1$ and $e_2$, it follows that $||g||_{T'}<||g||_T$.  This contradicts the fact that $g\in \mathbf U_T$ and that by assumption $||g||_T=||g||_{T'}$.
Hence $f$ is injective as claimed.

Thus the map $f$ is injective and is isometric on every edge of $T$. Therefore $f:T\to T'$ is an isometric embedding. Since the actions of $F_N$ on $T$ and $T'$ are minimal, it follows that $f(T)=T'$, so that $f:T\to T'$ is an $F_N$-equivariant isometry. Hence $T=T'$ in $\cvn$, as required.
 
\end{proof}
			
\section{Geodesic currents}\label{sect:currents}

\subsection{Basic facts.}

We will only state a few basic facts and definitions about currents, and refer the reader to~\cite{Martin,Ka1,Ka2,Ka3,KL1,KL2,KL3} for more detailed background information regarding geodesic currents.

For the free group $\FN$ define its ``double boundary'' $\partial^2 \FN$ as
\[
\partial^2 \FN:=\partial F_N \times \partial F_N - diag
=\{(\xi,\zeta)\in \partial F_N \times \partial F_N : \xi\ne \zeta\}. 
\]
The space $\partial^2 \FN$ comes equipped with a natural topology,
inherited from  and $\partial F_N \times \partial F_N$
natural translation action of $F_N$ by homeomorphisms. There is also a
natural ``flip'' map $\partial^2 \FN \to \partial^2 \FN$,
$(\xi,\zeta)\mapsto (\zeta,\xi)$, interchanging the two coordinates on
$\partial^2 \FN$.

Recall that a \emph{geodesic current} on $F_N$ is a positive Radon
measure $\mu$ on $\partial^2 F_N$ which is $F_N$-invariant and
flip-invariant. Here the ``flip'' map $\partial^2 F_N\to \partial^2
F_N$, $(X,Y)\mapsto (Y,X)$ interchanges the two
coordinates of $\partial^2 F_N$. The space $\Curr(F_N)$ of all geodesic
currents on $F_N$ comes equipped with a natural weak-* topology and a
natural left $\Aut(F_N)$-action by $\R_{\ge 0}$-linear homeomorphisms. The group of inner automorphisms of $F_N$ is contained in the kernel of this action, and hence the $\Aut(F_N)$ action naturally factors through to the action of $\Out(F_N)$ on $\Curr(F_N)$.

Every
nontrivial element $g\in F_N$ defines a \emph{counting} current
$\eta_g\in \Curr(F_N)$, which turns out to depend only on the conjugacy class $[g]$ of $g$ in $F_N$. 
Although the explicit definition of $\eta_g$ is not directly relevant for this paper, we briefly recall one of the equivalent definitions of $\eta_g$ here. Suppose first that $g\in F_N$ is a nontrivial element which is not a proper power in $F_N$. There are two well-defined distinct "poles" $g^\infty, g^{-\infty}\in \partial F_N$ where
\[
g^\infty =\lim_{n\to\infty} g^n, \quad g^{-\infty} =\lim_{n\to\infty} g^{-n}
\]
where the convergence is understood in the sense of the standard hyperbolic compactification $F_N\cup \partial F_N$ of $F_N$. Then $(g^{-\infty}, g^{\infty})\in \partial^2 F_N$.
Let $[g]$ denote the conjugacy class of $g$ in $F_N$.
Then
\[
\eta_g:=\sum_{h\in [g]\cup [g^{-1}]}\delta_{(h^{-\infty}, h^{\infty})}.
\]
Now if $g\in F_N$ is an arbitrary nontrivial element, $g$ can be uniquely written as $g=g_0^m$ where $m\ge 1$ and $g_0\in F_N$ is not a proper power. Then $\eta_g$ is defined as $\eta_g:=m\eta_{g_0}$.
We summarize the following basic facts about counting currents (see~\cite{Ka2}):

\begin{prop} Let $N\ge 2$ and let $g\in F_N$, $g\ne 1$.
Then
\begin{enumerate}
\item We have $\eta_g=\eta_{g^{-1}}$.
\item For every $n\in \mathbb Z$, $n\ne 1$, we have $\eta_{g^n}=n\eta_g$.
\item For every $h\in [g]$ we have $\eta_g=\eta_h$. Thus $\eta_g$ depends only on the conjugacy class of $g$, so we also use the notation $\eta_{[g]}:=\eta_g$.
\item For every $\phi\in \Out(F_N)$  we have
$\phi\eta_{[g]}=\eta_{\phi([g])}$.
\end{enumerate}

\end{prop}

The scalar multiples of counting currents are
called \emph{rational} currents. An important basic fact states that
the set $\mathcal R_N:=\{r\eta_g | r\ge 0, g\ne 1, g\in F_N\}$ of all rational currents is
dense in $\Curr(F_N)$. 
The space $\Curr(F_N)$ has a natural projectivization $\mathbb
PCurr(F_N)$ consisting of all equivalence classes $[\mu]$ where
$\mu\in Curr(F_N)$, $\mu\ne 0$. Here two currents $\mu_1,\mu_2$ are
equivalent if there exists $r>0$ such that $\mu_2=r\mu_1$.  The equivalence class $[\mu]$ is also called the \emph{projective class} of $\mu$.
The space
$\mathbb PCurr(F_N)$ is compact and it inherits a left action of
$\Out(F_N)$ by homeomorphisms. 
A basic fact about geodesic currents states (see~\cite{Ka2}):
\begin{prop}
Let $N\ge 2$ and let $g,h\in F_N$ be nontrivial elements. Then $[\eta_g]=[\eta_h]$ in $\PCN$ if and only if there exist $u\in F_N$, $m,n\in\Z$ such that $[g]=[u^m]$ and $[h]=[u^n]$.
\end{prop}

\subsection{The geometric intersection form.}

A key object connecting the Outer space and the space of geodesic currents is the so-called \emph{geometric intersection form}, constructed in~\cite{KL2}:

\begin{prop}\label{prop:gif}
Let $N\ge 2$.
There exists a unique continuous map 
\[
\langle \, \cdot \, , \, \cdot \, \rangle: \cvnbar\times Curr(F_N)\to \mathbb R_{\ge 0}
\]
with the following properties:
\begin{enumerate}
\item $\langle T, c_1\mu_1+c_2\mu_2\rangle=c_1\langle T, \mu_1\rangle+c_2\langle T,\mu_2\rangle$ for any $T\in \cvnbar$, $\mu_1,\mu_2\in Curr(F_N)$, $c_1,c_2\ge 0$. 
\item $\langle c T, \mu\rangle=c\langle T, \mu\rangle$ for any $T\in \cvnbar$, $\mu\in Curr(F_N)$ and $c\ge 0$.
\item $\langle \phi T, \phi\mu\rangle= \langle T,\mu\rangle$ for any $T\in \cvnbar$, $\mu\in Curr(F_N)$ and $\phi\in \Out(F_N)$.
\item $\langle T, \eta_g\rangle=||g||_T$ for every $T\in \cvnbar$ and $g\in F_N, g\ne 1$.
\end{enumerate}
\end{prop}
The value $\langle T, \mu\rangle$ is called the \emph{geometric intersection number} of $T\in \cvnbar$ and $\mu\in Curr(F_N)$.

\subsection{The minimal set}

\begin{defn}[Minimal set $\mathbb M_N$]
Let $N\ge 2$. Denote by $\mathbb M_N$ the closure in $\mathbb PCurr(F_N)$ of the set
\[
\{[\eta_g]: g\in F_N \text{ is primitive.}\}
\]
\end{defn}
It is easy to see that $\mathbb M_N\subseteq \mathbb PCurr(F_N)$ is a closed $\Out(F_N)$-equivariant subset. It turns out~\cite{KL1} that for $N\ge 3$ this is the minimal such subset:

\begin{prop}\label{prop:minset}
Let $N\ge 2$. Then:

\begin{enumerate}
\item~\cite{Martin} For every element $[\mu]\in \mathbb M_N$ the subset
  $\Out(F_N)[\mu]$ is dense in $\mathbb M_N$.
\item~\cite{KL1} Let $N\ge 3$. Then $\mathbb M_N\subseteq \mathbb PCurr(F_N)$ is the unique minimal closed $\Out(F_N)$-equivariant nonempty subset. This means that whenever $Z\subseteq  \mathbb PCurr(F_N)$ is a closed $\Out(F_N)$-equivariant nonempty subset then $\mathbb M_N\subseteq Z$. 

\end{enumerate}
\end{prop}

For $N\ge 3$ part (1) of Proposition~\ref{prop:minset} follows directly from part (2). For $N=2$ part (1) of Proposition~\ref{prop:minset} follows from the results of Reiner Martin~\cite{Martin} who showed that $\mathbb M_2$ is homeomorphic to the circle and that the action of $\Out(F_2)=GL(2,\Z)$ on $\mathbb M_2$ can be identified with the standard action of $GL(2,\Z)$ on $\mathbb S^1=\partial H^2$. 

The conclusion of part (2) of Proposition~\ref{prop:minset} is false for $N=2$. Indeed, for $F_2=F(a,b)$, it is easy to see that $\eta_{g}$ for $g=[a,b]$ is a fixed point  for the action of $\Out(F_2)$ on $\Curr(F_2)$ and hence $[\eta_g]$ is fixed by the action of  $\Out(F_2)$ on $\mathbb PCurr(F_2)$.  However, we will see later that a weaker version  of Proposition~\ref{prop:minset}  is true for $N=2$ and that version will be sufficient for our purposes.

\subsection{Consequences of spectral rigidity of the set of primitive elements}
Recall that $\mathcal P_N$ denotes the set of all primitive elements in $F_N$.

A key tool in proving our main results is the following:

\begin{prop}\label{prop:m}
Let $N\ge 2$, let $H\le Aut(F_N)$ and let $g\in F_N, g\ne 1$.
Suppose that  the closure of the set $H[\eta_g]$ in
$\PCN$ contains the set $\mathbb M_N$. 

Then \[
Hg=\{\phi(g)| \phi\in H\}  \subseteq F_N
\]
is a spectrally rigid subset of $F_N$.
\end{prop}

\begin{proof}
Let $T,T'\in \cvn$ be such that $||\phi(g)||_T=||\phi(g)||_{T'}$ for
every $\phi\in H$. We need to show that $T=T'$ in $\cvn$.

Let $Z$ be the closure of the set $H[\eta_g]$ in
$\PCN$.  By assumption we have $\mathbb
M_N\subseteq Z$. Hence for every primitive element $a\in F_N$ there exists a sequence $\phi_n\in H$
and a sequence $c_n\ge 0$ such that
\[
\lim_{n\to\infty} c_n\phi_n\eta_g=\lim_{n\to\infty} c_n\eta_{\phi_n(g)}=\eta_a \qquad \text{ in } Curr(F_N).
\] 
Proposition~\ref{prop:gif} then implies that
\begin{gather*}
||a||_T=\langle T, \eta_a\rangle= \langle T,
\lim_{n\to\infty} c_n\eta_{\phi_n(g)}\rangle =\lim_{n\to\infty}
c_n\langle T, \eta_{\phi_n(g)}\rangle=\\
\lim_{n\to\infty} c_n ||\phi_n(g)||_T=\lim_{n\to\infty} c_n ||\phi_n(g)||_{T'}=\lim_{n\to\infty}
c_n\langle T', \eta_{\phi_n(g)}\rangle=\\\langle T',
\lim_{n\to\infty} c_n\eta_{\phi_n(g)}\rangle=\langle T', \eta_a\rangle=||a||_{T'}.
\end{gather*}
Thus $||a||_T=||a||_{T'}$ for every primitive element $a\in F_N$.
Theorem~\ref{thm:prim} now implies that $T=T'$ in $\cvn$.
\end{proof}

\subsection{Stable and unstable currents}

Let $N\ge 2$. An element $\phi\in \Out(F_N)$ is called \emph{fully irreducible} or \emph{iwip} (for "irreducible with irreducible powers") if there is no integer $n\ge 1$ such that $\phi^n$ preserves the conjugacy class of a proper free factor of $F_N$.

For an element $\phi\in \Out(F_N)$ the conjugacy class $[g]$, where $g\in F_N$, $g\ne 1$, is \emph{periodic} if there exists $n\ge 1$ such that $\phi^n[g]=[g]$. An element $\phi\in \Out(F_N)$ is called \emph{atoroidal} if $\phi$ does not have any periodic conjugacy classes.

It is well-known that all non-atoroidal iwips in $\Out(F_N)$  come from homeomorphisms of compact surfaces with a single boundary component: 
\begin{prop}\cite{BH92}
Let $N\ge 2$ and $\phi\in \Out(F_N)$ be an iwip.  Then the following hold:

\begin{enumerate}

\item The automorphism $\phi$ is not atoroidal if and only if there exists an isomorphism $\alpha:F_N\to\pi_1(\Sigma)$, where $\Sigma$ is a connected compact surface with exactly one boundary component, such that $\phi$ is induced by a homeomorphism of $\Sigma$. 

\item Let $\phi\in \Out(F_N)$ be a non-atoroidal iwip, let $\Sigma$ be as in (1) and let $[h]$ be the conjugacy class in $F_N$ given by the boundary of $\Sigma$.

Then the only periodic conjugacy classes of $\phi$ in $F_N$ are those of the form $[h^m]$, where $m\in\Z, m\ne 0$. The conjugacy class $[h]$ is called the \emph{peripheral curve} of $\phi$. (Note that for an non-atoroidal iwip $\phi$ the peripheral curve $[h]$ is defined uniquely up to inversion. Namely, if $[g]$ is a periodic conjugacy class of $\phi$ such that $g\in F_N$ is not a proper power then $[g]=[h^{\pm 1}]$.)

\item Let $N=2$ and $F_2=F(a,b)$. Let $u=[a,b]\in F_2$.  Let $\phi\in \Out(F_2)$ be an iwip. Then $\phi$ is not atoroidal and $[u]$ is the peripheral curve of $\phi$.
\end{enumerate}
\end{prop}

For an element $\phi\in \Out(F_N)$ a current $\mu\in \Curr(F_N)$ is called an \emph{eigencurrent} of $\phi$ if $\mu\ne 0$ and $\phi\mu=\lambda\mu$ for some $\lambda\ge 0$. In that case the number $\lambda$ is called the associated \emph{eigenvalue} of $\mu$ for $\phi$. Thus for $\mu\in \Curr(F_N)$, $\mu\ne 0$ is an eigencurrent of $\phi$ if and only if $[\mu]\in \PCN$ is a fixed point of $\phi$. If $[\mu]\in\PCN$ is a fixed point of $\phi$, we also refer to the eigenvalue of $\mu$ for $\phi$ as the eigenvalue of $[\mu]$ for $\phi$, and we sometimes refer to $[\mu]$ as an eigencurrent of $\phi$.

Recall that $\mathbb M_N\subseteq \PCN$ is a closed
$\Out(F_N)$-invariant subset.
As proved by Reiner Martin~\cite{Martin}, if $\phi\in \Out(F_N)$ is an iwip, then $\phi$ has the "North-South" dynamics on the minimal set $\mathbb M_N\subseteq \PCN$ and, moreover, if $\phi$ is an atoroidal iwip, then  $\phi$ has the "North-South" dynamics on $\PCN$.  We only need the following weak version of Martin's result:

\begin{prop}[Stable and unstable eigencurrents of iwips]\cite{Martin}\label{prop:Martin}

Let $N\ge 2$ and let $\phi\in \Out(F_N)$ be an iwip.

Then the following holds:
\begin{enumerate}
\item The element $\phi$ has exactly two distinct fixed points in
  $\mathbb M_N$. One of these fixed points, called the \emph{stable eigencurrent} of $\phi$, and denoted by $[\mu_+]=[\mu_+(\phi)]\in \mathbb M_N$, has eigenvalue $>1$ for $\phi$, and the other fixed point, called the \emph{unstable eigencurrent} of $\phi$ and denoted by $[\mu_-]=[\mu_-(\phi)]\in \mathbb M_N$, has eigenvalue $<1$ for $\phi$.
Thus $\phi\mu_+=\lambda_+\mu_+$ for $\lambda_+>1$ and
$\phi\mu_-=\frac{1}{\lambda_-}\mu_-$ for $\lambda_->1$. 
\item If $\phi$ is both atoroidal and an iwip then for every
  $g\in F_N$, $g\ne 1$ we have
\[
\lim_{n\to\infty} \phi^n [\eta_g]=[\mu_+(\phi)], \qquad \lim_{n\to\infty} \phi^{-n} [\eta_g]=[\mu_-(\phi)].
\]
\item If $\phi$ is an iwip which is not atoroidal and if $[u]$ is the
  peripheral curve for $\phi$ then for every nontrivial $g\in F_N$ such that $g$
  is not conjugate to $u^m$, $m\in\Z$ we have
\[
\lim_{n\to\infty} \phi^n [\eta_g]=[\mu_+(\phi)], \qquad \lim_{n\to\infty} \phi^{-n} [\eta_g]=[\mu_-(\phi)].
\]
\item For any $\theta\in \Out(F_N)$ the element $\phi':=\theta
  \phi\theta^{-1}$ is again an iwip and $[\mu_+(\phi')]=\theta
  [\mu_+(\phi)]$, $[\mu_-(\phi')]=\theta
  [\mu_-(\phi)]$.
\end{enumerate}

\end{prop}

\section{Spectrally rigid automorphic orbits}

\subsection{Finding an iwip in an ample subgroup}

\begin{lem}\label{lem:iwip}
Let $N\ge 2$ and let $H\le \Out(F_N)$ be an infinite normal
subgroup. Then $H$ contains an iwip element $\phi$.
\end{lem}		
\begin{proof}
By a result of Handel and Mosher~\cite{HM} either $H$ contains an iwip or there exist a subgroup $H_1$ of finite index in $H$ and a proper free factor $B$ of $F_N$ such that every element of $H_1$ leaves the conjugacy class of $B$ invariant. 

Suppose that the latter case occurs, so that some subgroup $H_1$ of finite index in $H$ preserves the conjugacy class of a proper free factor $B$ of $F_N$. 

It is well known (see, for example,~\cite{BV1}) that $\Out(F_N)$ is virtually torsion-free and hence there are no infinite torsion subgroups in $\Out(F_N)$.  Thus if we show that every element of $H_1$ has finite order, this will imply that $H_1$ is finite, yielding a contradiction with the assumption that  $H$ is infinite and that  $H_1$ has finite index in $H$.

Note that if $\theta\in \Out(F_N)$ is arbitrary, then $\theta H_1 \theta^{-1}$ leaves the conjugacy class $[\theta(B)]$ invariant. Since $H$ is normal in $\Out(F_N)$,  it follows that $\theta H_1 \theta^{-1}$ has finite index in $\theta H \theta^{-1}=H$ and hence $H_\theta:=H_1\cap \theta H_1 \theta^{-1}$ has finite index in $H_1$. Thus every element of $H_1$ has a positive power belonging to $H_\theta$. Similarly, for any finite collection of elements $\theta_1, \dots, \theta_m\in \Out(F_N)$ the subgroup 
\[
H_{\theta_1, \dots, \theta_m}:=H_1\cap  \theta_1 H_1 \theta_1^{-1}\cap \dots \cap  \theta_m H_1 \theta_m^{-1}
\]
has finite index in $H_1$. Moreover, every element of $H_{\theta_1, \dots, \theta_m}$ leaves invariant each of the conjugacy classes $[B], [\theta_1(B)], \dots [\theta_m(B)]$ and every element of $H_1$ has a positive power that belongs to $H_{\theta_1, \dots, \theta_m}$.

Let $\psi\in H_1$ be arbitrary.
Choose a free basis $A=\{a_1,\dots, a_N\}$ of $F_N$ such that for some $1\le k\le N-1$ the set $\{a_1,\dots, a_k\}$ is a free basis of $B$.
Let $\Psi\in \Aut(F_N)$ be a lift of $\psi$ to $\Out(F_N)$. Since $\psi$ preserves the conjugacy class of $B$, for every $n\in\Z$ the cyclically reduced form of $\Psi^n(a_1)$ over $A$ does not involve $a_N^{\pm 1}$.  For each $m=2,\dots, N-1$ let $\theta_m$ be the automorphism of $F_N$ induced by the permutation of $A$ which interchanges $a_m$ and $a_N$ and leaves the other elements of $A$ fixed. Thus $\theta_m(B)$ is generated by a subset of $A$ that does not involve $a_m$. Let $n\ge 1$ be such that $\Psi^n$ belongs to $H_{\theta_2,\dots, \theta_{N-1}}$. Then the image $\Psi^n(a_1)$ is a freely reduced word over $A$ whose cyclically reduced form does not involve $a_2^{\pm 1}, \dots, a_N^{\pm 1} $, so that $\Psi(a_1)$ is conjugate to $a_1^{\pm 1}$ in $F_N$. By taking a larger finite collection of automorphisms $\theta$ of $F_N$ we can find a positive power $\Psi$ which sends every $a_i$ to a conjugate of $a_i$ in $F_N$, where $i=1,\dots, N$. 

Moreover, by considering the automorphisms $\theta$ of $F_N$ corresponding to all the elementary Nielsen transformations on $A$, we can find an even bigger positive power $\Psi^M$ of $\Psi$ with the property that $\Psi^M(u)$ is conjugate to $u$ in $F_N$ for every freely reduced word $u$ of length $\le 2$ over $A$.  It is well known and easy to see that  this implies that $\Psi^M$ is an inner automorphism of $F_N$, so that $\psi^M=1$ in $\Out(F_N)$.  Thus we have shown that every element of $H_1$ has finite order in $\Out(F_N)$, which, since $\Out(F_N)$ is virtually torsion-free, implies that $H_1$ is finite. However, this contradicts the assumption that $H$ is infinite and that $H_1$ has finite index in $H$.

Thus $H$ contains an iwip element, as required.
  
\end{proof}		

\begin{rem}
Note that for $N=2$ every automorphism of $F_2=F(a,b)$ has a periodic conjugacy class - namely the conjugacy class of the commutator $[a,b]$. Thus every iwip $\phi\in \Out(F_2)$ is toroidal. 

For $N\ge 3$ it is almost certainly the case that every nontrivial normal subgroup $H$ of $\Out(F_N)$ contains an atoroidal iwip (and not just an iwip, as Lemma~\ref{lem:iwip} proves). Such a strengthening of  Lemma~\ref{lem:iwip} would eliminate the need to consider Case~2 in the proof of Theorem~\ref{thm:mainN3} below. However, the presently available tools do not appear to be sufficient for establishing the existence of an atoroidal iwip in $H$.  Let $\phi\in H$ be an iwip whose existence is provided by Lemma~\ref{lem:iwip}. Suppose that $\phi$ is toroidal. We can use then ping-pong considerations for the action of $\Out(F_N)$ on the free factor complex $FF_N$ (now known to be Gromov-hyperbolic by a result of Bestvina and Feighn~\cite{BF11})) and choose a conjugate $\psi=\theta \phi \theta^{-1}\in H$ such that for all sufficiently large $n\ge 1$ the subgroup $\Gamma=\langle \phi^n, \psi^n\rangle\le H$ is free of rank two and that every nontrivial element of $\Gamma$ is again an iwip. One would then like to argue that (with the appropriate choices of $\theta$ and $n$) the element $\phi^n\psi^n\in H$ is an atoroidal iwip. The difficulty in proving this statement is that toroidal iwips do not act with "north-south" dynamics on $\PCN$.  Instead, a toroidal iwip $\beta$ has three distinct fixed points in $\PCN$: the stable (expanding) current $[\mu_+(\beta)]$,  the unstable (contracting) current $[\mu_-(\beta)]$ and the current $[\eta_u]$ (where $u$ is the peripheral curve for $\beta$) such that $\eta_u$ is fixed by $\beta$ in the non-projective sense. To prove that the iwip $\alpha=\phi^n\psi^n$ is atoroidal we need to show that there does not exist a nontrivial current in $Curr(F_N)$ that is fixed by $\alpha$. Establishing this fact requires first proving the following generalized north-south dynamical property for the action of any toroidal iwip  $\beta\in \Out(F_N)$ for $N\ge 3$: for any neighborhood $U$ of $[\mu_+(\beta)]$ in $\PCN$ and for any open set $V$ in $\PCN$ containing the segment between $[\mu_-(\beta)]$ and $[\eta_u]$ (where $u$ is the peripheral curve for $\beta$) there exists $M\ge 1$ such that for all $m\ge M$ we have $\beta^m(\PCN\setminus V) \subseteq U$.  This statement, which is most likely true, is not yet proved in the literature and establishing it requires more delicate arguments than those used by Reiner Martin~\cite{Martin} in the proof of "north-south" dynamics for the action of atoroidal iwips on $\PCN$.

\end{rem}

\subsection{The proof of Theorem~\ref{thm:A} for the case $N \geq 3$.}

\begin{lem}\label{lem:nc}
Let $N\ge 3$ and let $H\le \Out(F_N)$ be an infinite normal
subgroup. Then for every nontrivial element $g\in F_N$ there exists
$\psi\in H$ such that $\psi(g)$ is not conjugate to $g^{\pm 1}$
in $F_N$.
\end{lem}
\begin{proof}
Let $g\in F_N$, $g\ne 1$. Suppose that for every $\psi\in H$ 
$\psi(g)$ is conjugate to $g$ or $g^{-1}$ in $F_N$. Thus $H[\eta_g]
=[\eta_g]$. Hence for every $\theta\in \Out(F_N)$ the subgroup $\theta
H\theta^{-1}$ fixes $\theta[\eta_g]=[\eta_{\theta(g)}]$ in
$\PCN$. Since $H$ is normal in $\Out(F_N)$, it follows that $H$ fixes
$\theta[\eta_g]=[\eta_{\theta(g)}]$ for every $\theta\in \Out(F_N)$. 

Let $\varphi\in \Out(F_N)$ be any atoroidal iwip. Then $H$ fixes
$\varphi^n[\eta_g]$ for every $n\ge 1$. Hence, by
Proposition~\ref{prop:Martin}, we have
$H[\mu_+(\varphi)]=[\mu_+(\varphi)]$. By a result of~\cite{KL6}, for any atoroidal iwip $\varphi$
the stabilizer of $[\mu_+(\varphi)]$ in $\Out(F_N)$ is virtually cyclic
and contains $\langle \varphi \rangle$ as a subgroup of finite
index. Now choose two atoroidal iwips $\varphi_1, \varphi_2\in
\Out(F_N)$ so that $\langle \varphi_1\rangle\cap \langle
\varphi_1\rangle=\{1\}$. Then, by the above argument, the intersection of the stabilizers of $[\mu_+(\varphi_1)]$ and $[\mu_+(\varphi_2)]$ must be finite and must contain $H$, which contradicts the assumption
that $H$ is an infinite subgroup of $\Out(F_N)$.

\end{proof}

The following statement implies  Theorem~\ref{thm:A} from
the Introduction for the case $N\ge 3$:

\begin{thm}\label{thm:mainN3}
Let $N\ge 3$.  Let $H\le \Aut(F_N)$ be an ample subgroup.

Then for every $g\in F_N$, $g\ne 1$, the automorphic
orbit 
\[
Hg \subseteq F_N
\]
is a spectrally rigid subset of $F_N$.
\end{thm}
		
\begin{proof}

Denote by $\overline H\le \Out(F_N)$ the image of $H$ in
$\Out(F_N)$. Thus $\overline H$ is contains an infinite normal subgroup of
$\Out(F_N)$, and, without loss of generality we may assume that $\overline H$ is an infinite normal subgroup of
$\Out(F_N)$.

Let $g\in F_N$, $g\ne 1$ be arbitrary.
Let $Z$ be the closure of the set $H [\eta_g]$ in
$\PCN$. Proposition~\ref{prop:m} implies that in order to establish
that $Hg$ is a spectrally rigid
subset of $F_N$ it suffices to establish the following:

{\bf Claim.} We have $\mathbb M_N\subseteq Z$.

Lemma~\ref{lem:iwip} implies that $\overline H$ contains some iwip element
$\phi$. 

{\bf Case 1.} Suppose first that the iwip $\phi$ is atoroidal.

By Proposition~\ref{prop:Martin} we have
\[
\lim_{n\to\infty} \phi^n [\eta_g]=[\mu_+(\phi)].
\]
Therefore $[\mu_+(\phi)]\in Z$. Recall also that
$[\mu_+(\phi)]\in \mathbb M_N$.

Let $\theta\in \Out(F_N)$ be arbitrary.
Since $\overline H$ is
normal in $\Out(F_N)$, it follows that $\theta \phi\theta^{-1}\in \overline H$. 
The element $\phi'=\theta \phi\theta^{-1}$ is again an iwip. 
Hence the same argument as above implies that
$\mu_+(\theta\phi\theta^{-1})\in Z$.

Since $\mu_+(\theta\phi\theta^{-1})=\theta [\mu_+(\phi)]$, we see that
$\theta [\mu_+(\phi)]\in Z$ for every $\theta\in \Out(F_N)$.
By Proposition~\ref{prop:minset} the subset $\Out(F_N)
[\mu_+(\phi)]\subseteq \mathbb M_N$ is dense in $\mathbb M_N$. Since
$Z$ is closed, this implies that $\mathbb M_N\subseteq Z$, as claimed.

{\bf Case 2.} Suppose that the iwip $\phi$ is not atoroidal. Let $[u]$
be the peripheral curve of $[\phi]$.

Since $Hg=H\Psi(g)$ for every $\Psi\in H$, Lemma~\ref{lem:nc} implies
that, after possibly replacing $g$ by $\Psi(g)$ for some $\Psi\in H$,
we may assume that $[\eta_g]\ne [\eta_u]$. Then   
Proposition~\ref{prop:Martin} implies that
\[
\lim_{n\to\infty} \phi^n [\eta_g]=[\mu_+(\phi)].
\]
Thus again we see that $[\mu_+(\phi)]\in Z$.

Let $\theta\in\Out(F_N)$ be arbitrary and let
$\phi'=\theta\phi\theta^{-1}$. Thus
$[\mu_+(\phi')]=\theta[\mu_+(\phi)]$ and $\theta(u)$ is the peripheral
curve of $\phi'$.

If $[\eta_g]\ne [\eta_{\theta(u)}]$ then again Proposition~\ref{prop:Martin} implies that
$\lim_{n\to\infty} (\phi')^n
[\eta_g]=[\mu_+(\phi')]=\theta[\mu_+(\phi)]$, so that
$[\mu_+(\phi')]=\theta[\mu_+(\phi)]\in Z$.

Suppose now that $[\eta_g]=[\eta_{\theta(u)}]$. Then
Lemma~\ref{lem:nc} again implies that there is some $\psi\in \overline H$ such
that $\psi[\eta_g]=[\eta_{\psi(g)}]\ne [\eta_{\theta(u)}]$. Note that
  $\overline H[g]=\overline H[\psi(g)]$ since $\psi\in \overline H$.
Proposition~\ref{prop:Martin} now implies that 
\[
\lim_{n\to\infty} (\phi')^n [\eta_{\psi(g)}]=[\mu_+(\phi')].
\]
Thus again $[\mu_+(\phi')]=\theta[\mu_+(\phi)]\in Z$.

We have shown that $\theta[\mu_+(\phi)]\in Z$ for every
$\theta\in\Out(F_N)$. Since $\Out(F_N)[\mu_+(\phi)]$ is dense in
$\mathbb M_N$ and the set $Z$ is closed, it follows that $\mathbb
M_N\subseteq Z$.

Thus the Claim is verified, which, as noted above, via
Proposition~\ref{prop:m} now implies that $Hg$ is a spectrally rigid subset of $F_N$.

\end{proof}

\subsection{The case $N=2$}
		
Fix a free basis $A:=\{a,b\}$ of $F_2$, so that $F_2=F(a,b)$.

We need the following weaker version of
Proposition~\ref{prop:minset} for the case $N=2$:

\begin{prop}\label{prop:minset2} Let $g\in F_2, g\ne 1$ be such that $g$ is not conjugate
  to any integer power of $[a,b]$. Let $H\le \Out(F_2)$ be an infinite
  normal subgroup. Let $Z$ be any closed
  $H$-invariant subset of $\mathbb PCurr(F_2)$ such that
  $[\eta_g]\in Z$. Then $\mathbb M_2\subseteq Z$.
\end{prop}
\begin{proof}
Let $Z$ be the closure of the set $H[\eta_g]=\{[\eta_{\phi(g)}]: \phi\in
H\}$ in $\mathbb PCurr(F_2)$.

Lemma~\ref{lem:iwip} implies that there exists an iwip element
$\phi\in H$. Recall that in this case $u=[a,b]$ is the peripheral
curve for $\phi$ and, moreover, powers of $[a,b]$ are the only
periodic conjugacy classes for $\phi$ in $F_2$.

Since by assumption $g$ is not conjugate
  to any integer power of $[a,b]$, it follows from
  Proposition~\ref{prop:Martin} that
\[
\lim_{n\to\infty} \phi^n[\eta_g]=\lim_{n\to\infty}[\eta_{\phi^n(g)}]=[\mu_+(\phi)].
\]
Therefore $[\mu_+(\phi)]\in Z$. For every $\theta\in \Out(F_2)$
$\phi'=\theta\phi\theta^{-1}$ is again an iwip with
$[\mu_+(\phi')]=\theta[\mu_+(\phi)]$. The subgroup $H$ is normal in $\Out(F_N)$,  and therefore
$\phi'\in H$.  Now Proposition~\ref{prop:Martin} again implies that
$[\mu_+(\phi')]\in Z$. Since $[\mu_+(\phi')]=\theta[\mu_+(\phi)]$, we
have shown that $\theta[\mu_+(\phi)]\in Z$ for every $\theta\in
\Out(F_2)$. By Proposition~\ref{prop:minset},  every $\Out(F_2)$-orbit of a point of $\mathbb M_2$
is dense in $\mathbb M_2$,  and therefore $\mathbb M_2\subseteq Z$,
as required.
\end{proof}

\begin{rem}\label{rem:F2}
Suppose $g$ is conjugate in $F_2=F(a,b)$ to $[a,b]^k$ for some $k\in \Zb$. Then $\Aut(F_2)g\subseteq F_2$ is not spectrally rigid; hence any subset of $\Aut(F_2)g$ is not spectrally rigid. 

Indeed, let $T_A$ and $T_B$ be the Cayley graphs
of $F(a,b)$ with respect to $A=\{a,b\}$ and $B=\{a,ab\}$ accordingly.
As usual, we give all edges of $T_A, T_B$ length 1. Then $T_A,T_B\in
{\rm cv}_2$ and $T_A\ne T_B$ in ${\rm cv}_2$.  

It is well-known that for any $\phi\in \Aut(F(a,b))$, the element $\phi([a,b])$ is
conjugate to $[a,b]^{\pm 1}$ in $F(a,b)$.  It follows that for any free bases $A$ and $B$ of $F_2=F(a,b)$, and for any $\phi\in \Aut(F_2)$ we have $||\phi(g)||_A=||\phi(g)||_B=4|k|$.
The Cayley graphs $T_A$ and $T_B$ of $F_2$ with respect to $A$ and $B$ respectively are both points in ${\rm cv}_2$. Thus we see that for every $\phi\in \Aut(F_2)$ we have 
$||\phi(g)||_{T_A}=||\phi(g)||_{T_B}=4|k|$. Since we chose $A$ and $B$
so that $T_A\ne T_B$ in ${\rm cv}_2$, this shows that the orbit
$\Aut(F(a,b))g$ is not spectrally rigid in $F(a,b)$.  
\end{rem}

It turns out that Remark~\ref{rem:F2} provides the only obstruction to
extending Theorem~\ref{thm:mainN3} to the case $N=2$, and we obtain
the conclusion of Theorem~\ref{thm:A} from the Introduction for $N=2$:

\begin{thm}\label{thm:mainN2}
Let $F_2=F(a,b)$ and let $g\in F_2, g\ne 1$ be such that $g$ is not
conjugate to a power of $[a,b]$ in $F(a,b)$.

Let $H\le \Aut(F_2)$ be an ample subgroup.
Then $Hg$ is a 
spectrally rigid subset of $F(a,b)$.
\end{thm}

\begin{proof}
Suppose that $g\in F(a,b)$, $g\ne 1$ is such that $g$ is is not
conjugate to a power of $[a,b]$ in $F(a,b)$. Let $Z$ be the closure in
$\mathbb PCurr(F_N)$ of the set
$H[\eta_g]$. Proposition~\ref{prop:minset2} implies that $\mathcal
M_2\subseteq Z$. Proposition~\ref{prop:m} now implies that $Hg$ is a spectrally rigid
subset of $F_2$.

\end{proof}

\section{Open problems}\label{sect:op}

One can define a more restrictive notion of spectral rigidity than the
one considered in this paper. Namely, call a subset $S\subseteq F_N$
\emph{strongly spectrally rigid} if whenever $T, T'\in \cvnbar$ are
such that $||g||_T=||g||_{T'}$ for every $g\in S$ then $T=T'$ in
$\cvnbar$.

\begin{prob}\label{prob:strongrigidity}
For $N\ge 3$ is it true that every nontrivial $\Aut(F_N)$ orbit is
strongly spectrally rigid in $F_N$? Is it true that for $N\ge 3$ the set $\mathcal
P_N$ of all the primitive elements is strongly spectrally rigid in $F_N$?
\end{prob}

In~\cite{Ka4} it is proved that for the nonbacktracking simple random
walk on $F_N$ almost every trajectory of that walk gives a strongly
spectrally rigid subset of $F_N$. However, the proofs of
Theorem~\ref{thm:mainN3} and Theorem~\ref{thm:mainN2} in the present
paper do not imply strong spectral rigidity, primarily because our
proof of spectral rigidity of the set $\mathcal P_N$ in
Theorem~\ref{thm:prim} only works for the interior points of the Outer
space.  It was recently pointed out to us by Jing Tao that for $N=2$ the set $\mathcal P_2$ of all primitive elements in $F_2=F(a,b)$ is \emph{not} strongly spectrally rigid in $F_2$.
She provided the following example demonstrating this fact:

\begin{exmp}[Jing Tao's example]\label{exmp:Tao}
Consider a graph of groups $\mathbb A$ with the underlying graph $A$ consisting of two vertices $x_0$, $x_1$ and two distinct topological edges $e$ and $f$, each  with end-vertices $x_0$, $x_1$.  Thus $A$ is a topological circle subdivided in two edges.
Orient $e$ and $f$ so that both $e$ and $f$ have $x_0$ as the initial
vertex and $x_1$ as the terminal vertex.  Define the vertex and edge
groups as $A_{x_0}=A_{x_1}=A_{e}=\langle a\rangle$ and $A_f=\{1\}$,
with the boundary monomorphisms for the edge $e$ being the identity
map $\langle a\rangle \to \langle a\rangle$. Identify $b$ with the
loop $ef^{-1}$. This provides an identification between $F(a,b)$ and
$\pi_1(\mathbb A, x_0)$. Note that the graph of groups $\mathbb A$
gives a very small simplicial splitting of $F(a,b)$. Finally, for each
$t\in (0,1)$ give the edge $e$ length $t$ and the edge $f$ length
$1-t$. Thus for every $t\in (0,1)$ we get a point $T_t\in
\overline{\rm cv}_2$  given by the Bass-Serre tree of $\mathbb A$ with
the lifts of $e$ having length $t$ and the lifts of $f$ having length
$1-t$.  By construction, for any distinct $t,t'\in (0,1)$ we have
$T_t\ne T_{t'}$ in $\overline{\rm cv}_2$, and, moreover, the
projective classes of $T_t$ and $T_{t'}$ are also distinct. It is easy
to check that  for every  $t\in (0,1)$ and every cyclically reduced
word $w\in F_2$, where all $b$s occur with the same sign, the
translation length $||w||_{T_t}$ is equal to the absolute value of the
exponent sum on $b$ in $w$. Since every primitive element in $F(a,b)$
has cyclically reduced form where all occurrences of $a$ have the same
sign and all occurrences of $b$ have the same sign (see for
example~\cite{CMZ}), it follows that for any primitive element $w$ in
$F(a,b)$ the translation length $||w||_{T_t}$ is a constant function
of $t$.  Thus the restriction of $||.||_{T_t}$ to $\mathcal P_2$ is a
function that does not vary with $t$, which shows that the set
$\mathcal P_2$ is not strongly spectrally rigid. On the other hand,
for every $t\in (0,1)$ we have  $||[a,b]||_{T_t}=2t-2$, so that
knowing of the translation length of  $[a,b]$  is already sufficient
to distinguish $T_t$ from all the other trees in the family. As noted
by the referee, one can also obtain this example geometrically as
follows: Identify $F_2$ with the fundamental group $\pi_1(S)$ of a
torus with
a single boundary component, and put a hyperbolic metric on $S$
making the boundary component totally geodesic. Consider a measured
geodesic lamination $(L,\mu)$ on $S$ consisting of one closed geodesic (say the
meridian curve), of weight $t$, and a geodesic arc with two endpoints
on the boundary, of weight $1-t$. Then the tree dual to the lift of
this measured lamination to the universal cover $\widetilde S$ of $S$
is exactly the tree $T_t$ described above. 
\end{exmp}

It seems plausible, however, that the case of $F_2$ is special (because of the rather special nature of primitive elements in $F_2$)  and that for $N\ge 3$ the set $\mathcal P_N$ is strongly spectrally rigid in $F_N$. Some positive evidence in this direction is provided by the following observation.   Let $\phi\in \Out(F_N)$ (where $N\ge 3$) be an atoroidal iwip element and let  $T_\phi\in \cvn$ be the "stable tree" of an atoroidal iwip  $\phi\in \Out(F_N)$ (in particular $T_\phi \phi =\lambda T_\phi$, where $\lambda>1$ is the Perron-Frobenius eigenvalue of a train-track representative of $\phi$). We can prove that whenever $T'\in \cvn$ is such that the lengths functions of $T_\phi$ and of $T'$ agree on all primitive elements of $F_N$ then $T'=T_\phi$ in $\cvnbar$. Namely, in this case one can show that the Bestvina-Feighn-Handel "legal" lamination $L_{BFH}(\phi)$ of $\phi$ is contained in the dual algebraic lamination $L^2(T')$ of $T'$. This implies, for instance by the results of~\cite{KL3}, that  $[T_\phi]=[T']$ in $\CVNbar$, and it is then not hard to deduce that $T'=T_\phi$. We refer the reader to \cite{BFH97,KL3,CHL1,CHL2} for the background on dual algebraic laminations and on laminations associated to iwip automorphisms.

\begin{prob}
Let $N\ge 3$.  Does there exist a subset $S\subseteq F_N$ such that $S$ is spectrally rigid but not strongly spectrally rigid?
\end{prob}

\begin{prob}\label{prob:prim}
Let $N\ge 2$ and let $S\subseteq \mathcal P_N$ be an arbitrary subset of
the set $\mathcal P_N$ of all primitive elements in $F_N$. Since $S$
consists of primitive elements, this implies that  $\{[\eta_g]: g\in
S\}\subseteq \mathbb M_N$.

Is it true that $S$ is spectrally rigid in $F_N$ if and only if the
closure of $\{[\eta_g]: g\in S\}$ in $\PCN$ is equal to $\mathbb M_N$?
 
\end{prob}

It is easy to show, using the intersection form, that if the
closure of $\{[\eta_g]: g\in S\}$ in $\PCN$ is equal to $\mathbb M_N$
then $S$ is spectrally rigid. However, the converse implication
appears to be quite difficult. A recent result of 
Duchin, Leininger, and Rafi~\cite{DLR} establish a similar result to
that suggested in Problem~\ref{prob:prim} in the context of singular
flat metrics on surfaces.

We have seen in Theorem~\ref{thm:mainN3} that for $N\ge 3$ if $H\le \Aut(F_N)$ is a subgroup that projects to an infinite normal subgroup of $\Out(F_N)$ 
then the
$\Aut(F_N)$-orbit of every nontrivial element $g$ of $F_N$ is
spectrally rigid.  In an earlier version of this paper we conjectured that for every cyclic subgroup $H\le \Aut(F_N)$ and every $g\in F_N$ the orbit $Hg\subseteq F_N$ is not spectrally rigid in $F_N$. This conjecture was recently proved by Brian Ray~\cite{Ray}.
This fact and Theorem~\ref{thm:mainN3} suggests the following:

\begin{prob}\label{prob:H}
Let $N\ge 3$ and let $H\le \Aut(F_N)$ be an arbitrary subgroup. Is it
true that either for every nontrivial $g\in F_N$ the orbit $Hg\subseteq
F_N$ is spectrally rigid or that for every nontrivial $g\in F_N$ the orbit $Hg\subseteq
F_N$ is not spectrally rigid?
\end{prob}

A positive answer to the above question would mean that rigidity or
non-rigidity of the orbit $Hg$ (where $g\in F_N,g\ne 1$) depends only
on the subgroup $H\le \Aut(F_N)$ and not on the choice of a nontrivial element $g\in
F_N$.

Theorem~\ref{thm:C} motivates the following question:

\begin{prob}[Relatively strongly rigid finite sets]
Given $T\in \cvnbar$, does there exist a
finite subset $S\subseteq \mathcal P$ such that whenever
$T'\in\cvnbar$ is such that $||g||_T=||g||_{T'}$ for all $g\in S$ then
$T=T'$ in $\cvnbar$? What if we just require $S$ to be a finite subset
of $F_N$ (and not necessarily of $\mathcal P$)?
\end{prob}

As we noted in the introduction, by a result of
Cohen-Lustig-Steiner~\cite{CLS}, there does not exist a finite
spectrally rigid subset of $F_2$.  However, the argument in \cite{CLS}
involves looking at trees $T\in {\rm cv}_2$ with variable volume of
the quotient metric graph $T/F_2$. On the other hand, the
Smillie-Vogtmann construction~\cite{CV} for $N\ge 3$ only uses trees
with co-volume $1$, that is, points of $\CVN$.  In fact, for $N=2$ the
situation is quite different when restricting trees with quotient
graphs of volume $1$, that is, to ${\rm CV_2}$. Elaborating the arguments from~\cite{CV2} we can show
that for $F_2=F(a,b)$ the set $S_0:=\{a, b, ab, ab^{-1}, [a,b]\}$ is
"${\rm CV}_2$-rigid", that is, knowing the $||.||_T$-lengths, for an
arbitrary $T\in {\rm CV}_2$, of elements of $S_0$, uniquely determines
$T$.  This naturally leads to the following question:
\begin{prob}
Does there exist a finite ${\rm CV}_2$-rigid set of primitive elements in $F_2$?
\end{prob}

The notion of a spectrally rigid set naturally generalizes to the space of currents. Thus we say that a subset $S\subseteq Curr(F_N)$ is \emph{spectrally rigid} if whenever $T, T'\in \cvn$ are such that $\langle T, \mu\rangle=\langle T', \mu\rangle$ for every $\mu\in S$ then $T=T'$ in $\cvn$.  

\begin{prob}
Let $N\ge 2$. Does there exist a finite spectrally rigid set of currents $S\subseteq Curr(F_N)$?
\end{prob}

The original argument of Smillie-Vogtmann~\cite{SV} about non-existence of a finite spectrally rigid set of elements in $F_N$, and the above mentioned result of Brian Ray~\cite{Ray} about non-rigidity of orbits of cyclic subgroups of $\Aut(F_N)$ significantly rely on the fact that (in the language of currents) counting currents of elements of $F_N$ never have full support. Thus it is possible that there may indeed exist a finite spectrally rigid set of currents containing one or more current with full support. On the other hand, it also seems plausible that every finite subset $S\subseteq \mathbb M_N$ is not spectrally rigid.

\end{document}

\title{Corrigendum to: ``Spectral rigidity of automorphic orbits in free groups''}

\author{Mathieu Carette}
\address{\tt SST/IRMP, 
Chemin du Cyclotron 2, bte L7.01.01, 
1348 Louvain-la-Neuve,
Belgium}
\email{\tt mathieu.carette@uclouvain.be}

\author{Stefano Francaviglia}

\address{\tt Dipartimento di Matematica of the University of Bologna, 
Pizza Porta S. Donato 5, 40126 Bologna, Italy}
\email{\tt francavi@dm.unibo.it }

\author{Ilya Kapovich}
\address{\tt Department of Mathematics, University of Illinois at
 Urbana-Champaign, 1409 West Green Street, Urbana, IL 61801, USA
 \newline http://www.math.uiuc.edu/\~{}kapovich/} \email{\tt
  kapovich@math.uiuc.edu}

\author{Armando Martino}
\address{\tt School of Mathematics,
University of Southampton;
Highfield, Southampton,
SO17 1BJ, United Kingdom }
\email{\tt A.Martino@soton.ac.uk}

%\thanks{\today}

\begin{abstract}
Lemma 5.1 in our paper~\cite{CFKM} says that every infinite normal subgroup of $\Out(F_N)$ contains a fully irreducible element; this lemma was substantively used
in the proof of the main result, Theorem~A in \cite{CFKM}.
Our proof of Lemma 5.1 in \cite{CFKM} relied on a subgroup classification result of
Handel-Mosher~\cite{HM}, originally stated in \cite{HM} for arbitrary
subgroups $H\le \Out(F_N)$. It subsequently turned out (see p. 1 in \cite{HM1}) that the proof
of the Handel-Mosher theorem needs the assumption that $H$ be finitely
generated. Here we provide an alternative proof of Lemma~5.1 from
\cite{CFKM}, which uses the corrected version of the Handel-Mosher
theorem and relies on the 0-acylindricity of the action of $\Out(F_N)$
on the free factor complex (due to Bestvina-Mann-Reynolds).
\end{abstract}

\thanks{The first author is a Postdoctoral Researcher of the F.R.S.-FNRS (Belgium). The third author was supported by Collaboration Grant no. 279836 from the Simons Foundation}


\maketitle

%\tableofcontents

\section{Introduction}

The purpose of this note is to correct a gap in our
paper~\cite{CFKM}. The proof of the main result, Theorem~A, of~\cite{CFKM},
substantively relies on Theorem~1.1 of Handel and
Mosher~\cite{HM} about classification of subgroups of $\Out(F_N)$.
Originally Theorem~1.1 was stated in \cite{HM} for arbitrary subgroups
$H\le \Out(F_N)$, and we applied that statement in~\cite{CFKM}.  After
our paper~\cite{CFKM} was published, we were informed that the proof of Theorem~1.1 in \cite{HM} only goes through under
the additional assumption that $H\le\Out(F_N)$ be finitely generated;
see the footnote at p.1 of \cite{HM1}.

The specific use of Theorem~1.1 of \cite{HM} in \cite{CFKM} occurs in
the proof of Lemma~5.1 in \cite{CFKM}. This proof no longer works,
when the Handel-Mosher result is replaced by its finitely generated
version.  This situation has created a gap in the proof of Lemma~5.1 given
in \cite{CFKM}.

In this corrigendum we fix this gap and provide an alternative proof of Lemma~5.1. Thus Theorem~A in \cite{CFKM} and all the other results proved
there remain valid in their original form.  Lemma~5.1 in \cite{CFKM} stated the following:

\begin{prop}\label{prop:fix}
Let $N\ge 2$ and let $H\le \Out(F_N)$ be an infinite normal
subgroup. Then $H$ contains some fully irreducible element $\phi$.
\end{prop} 
Note that for $N\ge 3$ every nontrivial normal subgroup of
$\Out(F_N)$ is infinite, but $\Out(F_2)$ does possess a finite
nontrivial normal subgroup (namely the center of $\Out(F_2)$, which is
cyclic of order $2$). Recall also that an element $\phi\in\Out(F_N)$
is called \emph{fully irreducible} or \emph{iwip} (for ``irreducible
with irreducible powers'') if no positive power of $\phi$ preserves
the conjugacy class of a proper free factor of $F_N$.

The original formulation of Theorem~1.1 in \cite{HM} said that for an
arbitrary subgroup $H\le\Out(F_N)$, either $H$ contains a fully
irreducible element or $H$ has a subgroup of finite index $H_0$ such
that $H_0$ preserves the conjugacy class of some proper free factor of
$F_N$. As noted above, it turns out that the proof of
Theorem~1.1 in \cite{HM} only goes through under the additional
assumption that $H$ be finitely generated.

The new proof of Lemma~5.1 of \cite{CFKM}, presented here,
is quite different from our original argument in \cite{CFKM}, although the proof
still relies on the corrected finitely generated version of the
Handel-Mosher subgroup classification theorem. Another key
ingredient in this new argument is the proof, due to Bestvina, Mann
and Reynolds, of 0-acylindricity of the $\Out(F_N)$ action on the
free factor complex $\mathcal{FF}_N$.

The proof of 0-acylindricity was communicated to us
by Bestvina and Reynolds. Since this proof does not appear
anywhere in the literature, we include it here for completeness; see
Proposition~\ref{prop:acyl} below.

We are grateful to Ric Wade for bringing to our attention the issue
with the original formulation of Theorem~1.1 of \cite{HM}. We are also
grateful to Mladen Bestvina and Patrick Reynolds for explaining to us
the argument for establishing 0-acylindricity of
$\mathcal{FF}_N$. Finally, we thank Martin Lustig for very helpful
discussions regarding the train track theory.

\section{$0$-Acylindricity}

We will use the terminology and notations from~\cite{CFKM}. In
particular, $\cvn$ denotes the (unprojectivized) Outer space, $\CVN$
denotes the projectivized Outer space, $\cvnbar$ denotes the closure
of $\cvn$ in the hyperbolic length function topology, and $\CVNbar$ denotes the projectvization of $\cvnbar$, so that $\CVNbar$ is the standard compactification of $\CVN$.

Following \cite{Rey,BR}, if $T\in\cvnbar$ and $A\le F_N$ is a proper free factor of $F_N$, we say that $A$ \emph{reduces $T$} if there exists an $A$-invariant subtree $T'$ of $T$ such that $A$ acts on $T'$ with dense orbits (the subtree $T'$ is allowed to consist of a single point). For $T\in \cvnbar$ denote by $\mathcal R(T)$ the set of all proper free factors of $F_N$ which reduce $T$.  Note that in many cases $\mathcal R(T)$ is the empty set.

\subsection{Two useful lemmas}

Lemma~\ref{lem:reduce} and Lemma~\ref{lem:reduce1} below are
related statements, and either one  can be used to complete the
proof of 0-acylindricity. The proof of Lemma~\ref{lem:reduce} relies
on the use of relative train tracks and is based on an argument of
Brian Mann as communicated to us by Bestvina. The proof of
Lemma~\ref{lem:reduce1} relies on the use of geodesic currents and
algebraic laminations, and is due to Patrick Reynolds. We include both
lemmas here since each of them is interesting in its own right.

\begin{lem}\label{lem:reduce}
Let $h\in \Out(F_N)$ be an element of infinite order and let $A$ be a
proper free factor of $F_N$ such that $h[A]=[A]$. Then there exist a positive power $\psi=h^m$ of $h$ and a nontrivial
free factor $A'$ of $A$ with the following property:

 Let $T_0$ in $\cvn$ be arbitrary and let $n_i\to\infty$ be a
 subsequence such that for some $T\in\cvnbar$ we have $\lim_{i\to\infty} [T_0]\psi^{n_i}=[T]$ in $\CVNbar$.
 Then $A'$ reduces $T$.
\end{lem}

\begin{proof}
Among all the $h$-periodic, up to conjugacy, proper free factors of
$A$, choose a free factor $A'$ of smallest rank. Replace $h$ by a
positive power $\psi=h^m$ so that $\psi[A']=[A']$.

Let $T_0, n_i$ and $T$ be as in the statement of the lemma.

Now, using the construction of Bestvina-Handel~\cite{BH92},  there exists a relative train-track $f:\Gamma\to\Gamma$ for $\psi$ with
marking $\alpha:F_N\to\pi_1(\Gamma)$ such that the lowest
stratum $\Gamma_0$ of $\Gamma$ corresponds to $A'$.

Let $T_1\in \cvn$ be the point in Outer space determined by $\Gamma$
(with, say, the simplicial metric on $\Gamma$).

Then, after possibly passing to a further subsequence of $n_i$, we may assume that $[T_1]\psi^{n_i}$ converges to some tree $[S]\in \CVNbar$ as
$i\to\infty$.  We will first show that $A'$ reduces $S$.

For some $c_i\ge 0$ we have $c_iT_1\psi^{n_i}\to_{n\to\infty} S$ in $\cvnbar$.
Since $\psi\in\Out(F_N)$ has infinite order, we have $\lim_{n\to\infty}
c_i=0$.

If $A'$ has rank $1$, then $A'=\langle a \rangle$ and $\psi([a])=[a^{\pm
  1}]$. Since $c_i$ converges to $0$, we then have
$||a||_S=\lim_{i\to\infty} ||a||_{c_iT_1\psi^{n_i}}=\lim_{i\to\infty} c_i
||\psi^{n_i}(a)||_{T_1}=\lim_{i\to\infty} c_i ||a||_{T_1}=0$.
Thus $A'$ fixes a point of $S$ and hence $A'$ reduces $S$.

Suppose now that $A'$ has rank at least $2$. Then $\psi$ is
exponentially growing and the restriction of $\psi$ to $A'$ (after a
conjugation) is a fully irreducible automorphism $\psi'$ of $A'$. Let $\lambda$ be the largest P.F. eigenvalue among all the strata of $\Gamma$.
Theorem~6.2 of \cite{L} then implies that there exists $k\ge 0$ such that we can choose $c_i=\lambda^{-n_i}n_i^{-k}$ for the rescaling constants, 
so that $\lambda^{-n_i}n_i^{-k}T_1 \psi^{n_i}$ converges to (a scalar multiple of) $S$ in $\cvnbar$.  Let
$\lambda'$ be the P.F. eigenvalue for the bottom stratum of $\Gamma$,
that is $\lambda'>1$ is the stretch factor of the fully irreducible
$\psi'\in \Out(A')$. Thus for every nontrivial $a'\in A'$ (except possibly for $a'$ conjugate to a power of the peripheral curve, if $\psi'$ is a  non-hyperbolic fully irreducible and is induced by a pseudo-Anosov of a compact surface with a single boundary component)  we see that $||\psi^n(a')||_{T_1}$ grows as $(\lambda')^n$.

If $\lambda'<\lambda$ then for each $a'\in A'$  $||\psi^n(a')||_{T_1}$ grows as at a smaller exponential rate than $\lambda^n$ and hence, by the above choice of the rescaling constants $c_i$, we get
 $||a'||_S=0$. Thus in this case $A'$ fixes a point of $S$, so that $A'$ reduces $S$.

Suppose now that $\lambda=\lambda'$.  

If $k>0$, then, since $c_i=\lambda^{-n_i}n_i^{-k}$, we see that for every $a'\in A'$ the sequence  $\lambda^{-n_i}n_i^{-k}||a'||_{T_1 \psi^{n_i}}=\lambda^{-n_i}n_i^{-k}||\psi^{n_i}(a)||_{T_1}$ is bounded above by $c' n_i^{-k}\to_{i\to\infty} 0$ (where $c'>0$ is some constant depending on $a'$) and hence $||a'||_S=0$.  Again, in this case $A'$ fixes a point of $S$ and thus reduces $S$.

Finally, suppose that $k=0$, so that $c_i=\lambda^{-n_i}$.
Then for every $a'\in A'$ we have $||a'||_S=\lim_{i\to\infty} \frac{||\psi^{n_i}(a') ||_{T_1}}{\lambda^{n_i}}=||a'||_{Y_+}$ where $[Y_+]$ in
the compactified Outer space  of $A'$ is the attracting tree of
$\psi'$. Therefore the minimal $A'$-invariant subtree $S_{A'}$ of $S$ is
$A'$-equivariantly isometric to $Y_+$. Since $\psi'$ is a fully
irreducible automorphism of $A'$, the action of $A'$ on $Y_+=S_{A'}$ has
dense orbits (see \cite{LL}, for example). Thus $A'$ reduces $S$.

Now we return to the original tree $T_0\in\cvn$ such that $[T_0]\psi^{n_i}\to_{i\to\infty} [T]\in\CVNbar$.
Since we also have $[T_1]\psi^{n_i}\to_{i\to\infty} [S]$, a general result of Sela~\cite{Sela} implies that  there exists a
bi-Lipschitz $F_N$-equivariant homeomorphism $f: S\to T$. 
The fact that $A'$ reduces $S$ now implies that $A'$ also reduces $T$.  Indeed,
let $Y$ be an $A$-invariant subtree of $S$ such that $A'$ acts on $Y$ with dense orbits. Then $f(Y)$ is a connected subset of $T$ and hence $f(Y)$ is a subtree of $T$.
Since $A'$ acts on $Y$ with dense orbits and $f$ is bi-Lipschitz and $F_N$-equivariant, it follows that $A'$ acts on $f(Y)$ with dense orbits as well. Thus $A'$ reduces $T$, as claimed.
\end{proof}

If a tree $T\in\cvnbar$ does not have dense $F_N$-orbits, then it is known (see \cite{CHL2,Rey,KL3}) that $T$ canonically decomposes as a "graph of actions".
In this case there exists an $F_N$-equivariant distance non-increasing map $f:T\to Y$ for some very small simplicial metric tree $Y\in\cvnbar$ such for every vertex $v$
  of $Y$ the stabilizer $Stab_{F_N}(v)$ acts with dense orbits on some $Stab_{F_N}(v)$-invariant subtree $T_v$ of $T$ (where $T_v$ may be a single point). Moreover, the tree $Y$ is obtained from $T$ be collapsing each $T_v$ to a point, where $v$ varies over all vertices of $Y$.  In this situation we will say that  $Y\in\cvn$ is the \emph{simplicial tree associated} to $T$. 
Note that if $T\in\cvn$ then $Y=T$.

\begin{lem}\label{lem:reduce1}
Let $A$ be a proper free factor of $F_N$ and let $h_n\in \Out(F_N)$ be an infinite sequence of distinct elements of $\Out(F_N)$ such that $h_n([A])=[A]$ for all $n\ge 1$.
Let $T_0\in\cvn$ and $T\in\cvnbar$ be such that $[T_0]h_n\to [T]$ in
$\CVNbar$ as $n\to\infty$. Then:
\begin{enumerate}
\item If $T$ has dense $F_N$-orbits then there exists a nontrivial free factor $A'$ of $A$ such that $A'$
 reduces $T$.

\item If $T$ does not have dense $F_N$-orbits and $Y\in\cvnbar$ is the associated simplicial tree, then some
  nontrivial free factor $A'$ of $A$ reduces $Y$.

\end{enumerate}

\end{lem}
\begin{proof}

There exist $c_n\ge 0$ such that $\lim_{n\to\infty} c_n T_0h_n=T$ in $\cvnbar$. Since the elements $h_n$ are distinct and the action of $\Out(F_N)$ on $\CVN$ is properly discontinuous, it follows that $T\in \cvnbar\setminus\cvn$ and that $\lim_{n\to\infty}c_n=0$.
For every $n\ge 1$ choose a representative $\beta_n\in \Aut(F_N)$ of the outer automorphism $h_n$ such that $\beta_n(A)=A$.

Choose a nontrivial element $a\in A$ and put $a_n=\beta_n^{-1}(a)\in A$.
Then 
\[
c_n||a_n||_{T_0h_n} =||a_n||_{c_nT_0h_n}=||a_n||_{c_nT_0\beta_n}=c_n||\beta_n(a_n)||_{T_0}=c_n||a||_{T_0}\to_{n\to\infty} 0.
\]

In $\PCN$ we have $\lim_{n\to\infty}
\frac{1}{||a_n||_{T_0}}\eta_{a_n}=\mu\ne 0$.  
Here, for a nontrivial element $g\in F_N$, $\eta_g\in \Curr(F_N)$ is the
``counting current'' associated to $g$~\cite{Ka}.

Since $T_0\in\cvn$, there exists $c>0$ such that $||g||_{T_0}\ge c>0$ for all $g\in F_N$.
Hence $||a_n||_{T_0}\ge c$ and therefore
\begin{gather*}
\langle T,\mu\rangle=\lim_{n\to\infty} \langle c_nT_0h_n,
\frac{1}{||a_n||_{T_0}}\eta_{a_n}\rangle=
\lim_{n\to\infty} \frac{c_n}{||a_n||_{T_0}} \langle T_0h_n, \eta_{a_n}\rangle=\\
\lim_{n\to\infty} \frac{c_n}{||a_n||_{T_0}} ||a_n||_{T_0h_n}=\lim_{n\to\infty} \frac{c_n}{||a_n||_{T_0}}||a||_{T_0} = 0.
\end{gather*}
Here $\langle\ ,\ \rangle:\cvnbar\times Curr(F_N)\to \mathbb R$ is the
continuous ``geometric intersection form'' constructed in \cite{KL2}.

Since $\langle T,\mu\rangle=0$, by the main result of \cite{KL3} we have $supp(\mu)\subseteq
L(T)$, where $supp(\mu)$ is the support of $\mu$ and $L(T)$ is the
dual algebraic lamination of $T$ (see \cite{CHL2} for background about dual algebraic laminations associated to elements of $\cvnbar$). Since $a_n\in A$ for all $n\ge 1$, the construction of $\mu$ implies that there is a leaf of $L(T)$ that is carried by $A$.

If $T$ has dense $F_N$-orbits, then by Corollary~6.7 of~\cite{Rey},
there exists a nontrivial free factor $A'$ of $A$ such that $A'$
reduces $T$, and part (1) of the lemma holds.

Suppose now that $T$ does not have dense $F_N$-orbits, and let $Y\in\cvnbar$  be the associated simplicial tree for $T$.

Then by Lemma~10.2 of \cite{KL3}, we have $L(T)\subseteq
L(Y)$.  Since the factor $A$ carries a leaf of $L(T)$, it follows
that $A$ also carries a leaf $\ell$ of $L(Y)$. The  description of the dual lamination of a very small simplicial tree,
given in Lemma~8.2 of \cite{KL3}, then implies that $A$
contains some nontrivial element acting elliptically on $Y$. Namely, consider a free basis $X$ of $F_N$ such that $X$ contains as a subset a free basis $X'$ of $A$.  Lemma~8.2 of \cite{KL3} now implies that for some vertex group $U$ of $T'$ and for the Stallings core subgroup graph $\Delta_U$~\cite{KM} (with oriented edges labelled by elements of $X^{\pm 1}$) representing the conjugacy class of $U$, there exists a bi-infinite reduced edge-path $\gamma$ in $\Delta_U$ corresponding to the leaf $\ell$ of $L(Y)$. The fact that $\ell$ is carried by $A$ means that all the edges of $\gamma$ are labelled by elements of $(X')^{\pm 1}$. Since $\gamma$ is an infinite reduced path, we can find an immersed circuit as a subpath of $\gamma$. Then the label $a'$ of this circuit is a nontrivial element of $A$ whose conjugate belongs to $U$ and thus $a'$ acts elliptically on $Y$.

If $A$ fixes a point of $Y$, then $A$ reduces $Y$ and we are
done. Otherwise consider the minimal $A$-invariant subtree $Y_A$ of
$Y$. Then the quotient graph of groups $\mathbb A=Y_A//A$ gives a
nontrivial very small splitting of $A$ with at least one nontrivial
vertex group. The general structural result (Lemma 4.1 of~\cite{BF93}) about very small
simplicial splittings of free groups then implies that there exists a nontrivial free factor $A'$ of $A$ such that
a conjugate of $A'$ is contained in some vertex group of $\mathbb
A$. Then $A'$ fixes a vertex of the tree $Y_A$ and therefore $A'$
reduces $Y$, as required.

\end{proof}

\subsection{0-acylindricity}

For $N\ge 2$ let $\mathcal{FF_N}$ be the free factor graph for
$F_N$. The vertices of $\mathcal{FF_N}$ are the conjugacy classes
$[A]$ of proper free factors $A$ of $F_N$. For $N\ge 3$ the adjacency
of vertices in $\mathcal F_N$ corresponds to containment: two distinct
vertices $[A]$ and $[B]$ are adjacent if there exist representatives
$A$ of $[A]$ and $B$ of $[B]$ such that $A\le B$ or $B\le A$. For
$N=2$ the definition of adjacency is somewhat different, see
\cite{BF11} for details. It is known that for $N\ge 2$ the graph
$\mathcal{FF}_N$ is connected and Gromov-hyperbolic~\cite{BF11,KRa,HH}. There is a
natural action of $\Out(F_N)$ on $\mathcal{FF}_N$ by simplicial isometries.

We denote the simplicial metric on $\mathcal{FF}_N$ by $d$.

\begin{prop}[0-acylindricity of the free factor complex]\label{prop:acyl}
There exists a constant $M\ge 1$ with the following property:

If $N\ge 2$ and if $A,B$ are proper free factors of $F_N$ such that $d([A],[B])>M$ then
\[
Stab_{\Out(F_N)}([A]) \cap Stab_{\Out(F_N)}([B]) 
\]
is finite.
\end{prop}

\begin{proof}
Corollary~5.3 of \cite{BR} implies that there exists a constant $C>0$ (independent of the rank $N$ of $F_N$)  such that if $T\in\cvnbar$ admits a reducing factor then the set $\mathcal R(T)$ of all reducing factors for $T$ has diameter $\le C$ in $\mathcal{FF}_N$.  Take $M=C+2$. Let $A,B$ be proper free factors of $F_N$ such that $d([A],[B])>M$.
Put $H:=Stab_{\Out(F_N)}([A]) \cap Stab_{\Out(F_N)}([B])$. We claim that $H$ is finite.

Indeed, suppose not, and $H$ is infinite. Since $\Out(F_N)$ is virtually torsion-free, it follows that there exists an element $h\in H$ of infinite order.
Let $T_0\in \cvn$ be arbitrary and let $T\in\cvnbar$ be such that
$\lim_{i\to\infty} [T_0h^{n_i}]=[T]$ for some subsequence $n_i\to\infty$.

Since $h[A]=[A]$ and $h[B]=[B]$, by Lemma~\ref{lem:reduce} or,
alternatively, by Lemma~\ref{lem:reduce1}, there exist a tree
$S\in\cvnbar$ and
nontrivial free factors $A'$ of $A$ and $B'$ of $B$ such that $A'$ and
$B'$ both reduce $S$. When appealing to Lemma~\ref{lem:reduce}, we can
take $S=T$. If appealing to Lemma~\ref{lem:reduce1}, we can take $S=T$
if $T$ has dense $F_N$ orbits and otherwise we can take $S=Y$ where
$Y\in\cvnbar$ is the simplicial tree associated to $T$.

Note that if $N=2$ then $A'=A$, $B'=B$, and the
factors $A,B$ are infinite cyclic.

Then $d([A'],[A])\le 1$, $d([B'],[B])\le 1$ and therefore $d([A'],[B'])> M-2=C$.
Thus the set of reducing factors for $S$ has diameter $>C$ in $\mathcal{FF}_N$, which contradicts the choice of $C$.
\end{proof}

\section{The proof of Proposition~\ref{prop:fix}}
We can now recover Lemma~5.1 of \cite{CFKM}, which is stated as
Proposition~\ref{prop:fix} above.

\begin{proof}[Proof of  Proposition~\ref{prop:fix}]
Suppose that $H\le Out(F_N)$ is an infinite normal subgroup but that $H$ does not contain a fully irreducible element.
Since $\Out(F_N)$ is virtually torsion-free and $H$ is infinite, it  follows that H contains an element $\phi$ of infinite order.

By assumption on $H$, $\phi$ is not fully irreducible and hence, after replacing $\phi$ by a positive power, $\phi$ fixes the conjugacy class $[A]$ of a proper free factor $A$ of $F_N$.

Now let $M\ge 1$ be the 0-acylindricity constant provided by Proposition~\ref{prop:acyl}. We choose a fully irreducible $\theta\in \Out(F_N)$ and look at the conjugates $\alpha_n= \theta^n \phi\theta^{-n}$.
Note that $\alpha_n$ fixes the conjugacy class $\theta^n[A]$.  Since $H$ is normal, we have $\alpha_n\in H$, and thus the subgroup $L_n=\langle \phi, \alpha_n\rangle$ is contained in $H$.
The subgroup $L_n$ of $\Out(F_N)$ is finitely generated so the (corrected)
Handel-Mosher subgroup classification theorem~\cite{HM,HM1} does apply to $L_n$.

Since we assumed that $H$ contains no fully irreducible elements, $L_n$ must contain a subgroup $K_n$ of finite index which preserves a vertex $[B_n]$ of the free factor complex. Hence some positive powers of $\phi$ and of $\alpha_n$ preserve $[B_n]$.

Since $\theta$ is fully irreducible, we have $d([A], \theta^n[A])\to\infty$  as $n\to\infty$ (see \cite{BF08,KL2}). We choose $n\ge 1$ big enough so that $d([A], \theta^n[A]) > 2M+1$.

Then either $d([A], [B_n])>M$ or $d([B_n], \theta^n[A])>M$.

In the first case we get that some positive power $\phi^i$ of $\phi$ fixes both the vertices $[A]$ and $[B_n]$ of $\mathcal{FF}_N$, so that $\phi^i$ belongs to the intersection of their stabilizers. This contradicts 0-acylindricity since $\phi$ has infinite order.

In the second case for some $i>0$ the element $\alpha_n^i=\theta^n\phi^i\theta^{-n}\in H$ fixes both $[B_n]$ and $\theta^n[A]$. Thus $\alpha_n^i$ belongs to the intersection of the stabilizers of $[B_n]$ and $\theta^n[A]$, which again contradicts 0-acylindricity, since $\alpha_n$ has infinite order.
\end{proof}

\begin{rem}
Let $\phi\in H$ be a fully irreducible element provided by
Proposition~\ref{prop:fix}. Theorem~8.5 of Dahmani, Guirardel and
Osin~\cite{DGO} then implies that for some $m\ge 1$ the normal closure
$U=ncl(\phi^m)$ of $\phi^m$ in $\Out(F_N)$ is free of infinite rank
and every nontrivial element of $U$ is fully irreducible. Since $H\le
\Out(F_N)$ is normal by assumption, we have $U\le H$.
\end{rem}